\documentclass{amsart}
\usepackage{amssymb}
\usepackage{amscd}
\usepackage{verbatim}
\usepackage{epsfig}

\newcommand\Deltat{\tilde\Delta}
\begin{document}

\newcommand\Mand{\ \text{and}\ }
\newcommand\Mwith{\ \text{with}\ }
\newcommand\Mfor{\ \text{for}\ }
\newcommand\Mst{\ \text{such that}\ }
\newcommand\Mor{\ \text{or}\ }
\newcommand\Mif{\ \text{if}\ }
\newcommand\Miff{\ \text{iff}\ }
\newcommand\Mthen{\ \text{then}\ }
\newcommand\nin{\notin}
\newcommand\identity{\operatorname{id}}
\newcommand\Id{\operatorname{Id}}
\newcommand\Real{\mathbb{R}}
\newcommand\RR{\mathbb{R}}
\newcommand\CC{\mathbb{C}}
\newcommand\olM{\overline{M}}
\newcommand\olX{\overline{X}}
\newcommand\pos{\Real^+}
\newcommand\Rnp{\Real\setminus\{0\}}
\newcommand\nzero{\setminus\{0\}}
\newcommand\Cx{\mathbb{C}}
\newcommand\Cxp{\Cx^+}
\newcommand\Cxm{\Cx^-}
\newcommand\Nat{\mathbb{N}}
\newcommand\halfNat{{\frac{1}{2}}\mathbb{N}}
\newcommand\intgr{\mathbb{Z}}
\newcommand\im{\operatorname{Im}}
\newcommand\re{\operatorname{Re}}
\newcommand\sign{\operatorname{sign}}
\newcommand\codim{\operatorname{codim}}
\newcommand\End{\operatorname{End}}
\newcommand\Ker{\operatorname{Ker}}
\newcommand\Hom{\operatorname{Hom}}
\newcommand\tr{\operatorname{tr}}
\newcommand\Tr{\operatorname{Tr}}
\newcommand\ideal{{\mathcal I}}
\newcommand\Span{\operatorname{span}}
\newcommand\image{\operatorname{image}}
\newcommand\Range{\operatorname{Ran}}
\newcommand\Graph{\operatorname{graph}}
\newcommand\slim{\operatornamewithlimits{s-lim}}
\newcommand\gll{\mathfrak{gl}}
\newcommand\sll{\mathfrak{sl}}
\newcommand\sol{\mathfrak{so}}
\newcommand\GL{\operatorname{G\ell}}
\newcommand\SL{\operatorname{SL}}
\newcommand\SO{\operatorname{SO}}
\newcommand\On{\operatorname{O}}
\newcommand\pa{\partial}
\newcommand\del{\partial}
\newcommand\Rn{\Real^n}
\newcommand\Rm{\Real^m}
\newcommand\RN{\Real^N}
\newcommand\RtN{\Real^{2N}}
\newcommand\RM{\Real^M}
\newcommand\HH{\mathbb{H}}
\newcommand\sphere{\mathbb{S}}
\newcommand\Sn{\sphere^{n-1}}
\newcommand\Sm{\sphere^{m-1}}
\newcommand\Snp{\sphere^n_+}
\newcommand\Smp{\sphere^m_+}
\newcommand\SN{\sphere^{N-1}}
\newcommand\SNp{\sphere^N_+}
\newcommand\circlep{\sphere^1_+}
\newcommand\Phom{P_{h}}
\newcommand\Shom{S_{h}}
\newcommand\distance{\operatorname{dist}}
\newcommand\cl{\operatorname{cl}}
\newcommand\interior{\operatorname{int}}
\newcommand\Fa{\operatorname{Fa}}
\newcommand\ff{\operatorname{ff}}
\newcommand\mf{\operatorname{mf}}
\newcommand\cf{\operatorname{cf}}
\newcommand\scf{\operatorname{sf}}
\newcommand\lf{\operatorname{lf}}
\newcommand\rf{\operatorname{rf}}
\newcommand\indfam{{\mathcal K}}
\newcommand\ev{{\lambda}}
\newcommand\fraka{{\mathfrak a}}
\newcommand\frakb{{\mathfrak b}}
\newcommand\frakw{{\mathfrak w}}
\newcommand\frakg{{\mathfrak g}}
\newcommand\frakk{{\mathfrak k}}
\newcommand\frakp{{\mathfrak p}}
\newcommand\calA{{\mathcal A}}
\newcommand\calB{{\mathcal B}}
\newcommand\calR{{\mathcal R}}
\newcommand\calS{{\mathcal S}}
\newcommand\calO{{\mathcal O}}
\newcommand\calJ{{\mathcal J}}
\newcommand\calM{{\mathcal M}}
\newcommand\calN{{\mathcal N}}
\newcommand\calX{{\mathcal X}}
\newcommand\calY{{\mathcal Y}}
\newcommand\calF{{\mathcal F}}
\newcommand\calG{{\mathcal G}}
\newcommand\calT{{\mathcal T}}
\newcommand\calC{{\mathcal C}}
\newcommand\calP{{\mathcal P}}
\newcommand\calU{{\mathcal U}}
\newcommand\calV{{\mathcal V}}
\newcommand\calCt{{\tilde {\mathcal C}}}
\newcommand\calCL{{\mathcal C}_{\text L}}
\newcommand\calCR{{\mathcal C}_{\text R}}
\newcommand\Cinf{{\mathcal C}^{\infty}}
\newcommand\dist{{\mathcal C}^{-\infty}}
\newcommand\dCinf{\dot{\mathcal C}^\infty}
\newcommand\ddist{\dot\dist}
\newcommand\Cj{{\mathcal C}^j}
\newcommand\Linf{L^{\infty}}
\newcommand\phg{{\text{phg}}}
\newcommand\bcon{{\mathcal A}}
\newcommand\bconc{{\mathcal A}_{\text{phg}}}
\newcommand\Sch{{\mathcal S}}
\newcommand\temp{\Sch^{\prime}}
\newcommand\Diff{\operatorname{Diff}}
\newcommand\Diffb{\operatorname{Diff}_{\text{b}}}
\newcommand\Diffc{\operatorname{Diff}_{\text{c}}}
\newcommand\Diffsc{\operatorname{Diff}_{\text{sc}}}
\newcommand\DiffI{\operatorname{Diff}_{\text{I}}}
\newcommand\DiffIq{\operatorname{Diff}_{\text{I},q}}
\newcommand\sing{\text{sing}}
\newcommand\reg{\text{reg}}
\newcommand\supp{\operatorname{supp}}
\newcommand\ssupp{\operatorname{sing\ supp}}
\newcommand\csupp{\operatorname{cone\ supp}}
\newcommand\esupp{\operatorname{ess\ supp}}
\newcommand\Fr{{\mathcal F}}
\newcommand\Frinv{\Fr^{-1}}
\newcommand\bop{{\mathcal B}}
\newcommand\spec{\operatorname{spec}}
\newcommand\pspec{\spec_{pp}}
\newcommand\cspec{\spec_{c}}
\newcommand\FIO{{\mathcal I}}
\newcommand\SP{\operatorname{RC}}
\newcommand\RC{\operatorname{RC}}
\newcommand\Symc{S_c}
\newcommand\Symca{S_c^{\alpha}}
\newcommand\Symczero{S_c^{0,...,0}}
\newcommand\sci{{}^{\text{sc}}}
\newcommand\sct{\sci T^*}
\newcommand\scdt{\sci \dot T^*}
\newcommand\dS{\dot S^*}
\newcommand\dT{\dot T^*}
\newcommand\dSreg{\dot\Sigma_{\text reg}}
\newcommand\scct{\sci\bar{T}^*}
\newcommand\Csc{C_{\text{sc}}}
\newcommand\SNpscd{(\SNp)^2_{\text{sc}}}
\newcommand\scdiag{\Delta_{\text{sc}}}
\newcommand\projscl{\pi^L_{\text{sc}}}
\newcommand\projscr{\pi^R_{\text{sc}}}
\newcommand\scHL{\sci H^{2,0}_{|\zeta|^2-\lambda^2}}
\newcommand\scHrg{\sci H^{2,0}_{\sqrt{g}}}
\newcommand\Hsc{H_{\text{sc}}}
\newcommand\WF{\operatorname{WF}}
\newcommand\WFp{\operatorname{WF^{\prime}}}
\newcommand\WFsc{\operatorname{WF}_{\text{sc}}}
\newcommand\WFscp{\operatorname{WF_{sc}^{\prime}}}
\newcommand\WFC{\operatorname{WF}_C}
\newcommand\WFCi{\operatorname{WF}_{C_i}}
\newcommand\elliptic{\operatorname{ell}}
\newcommand\Psop{\operatorname{\Psi}}
\newcommand\Psiscrs{\operatorname{\Psi_{sc}^{-2,\infty}}}
\newcommand\Psiscr{\operatorname{\Psi_{sc}^{-2,0}}}
\newcommand\Psiscrm{\operatorname{\Psi_{sc}^{0,2}}}
\newcommand\PsiscHam{\operatorname{\Psi_{sc}^{2,0}}}
\newcommand\Psisci{\operatorname{\Psi_{sc}^{*,*}}}
\newcommand\Psiscid{\operatorname{\Psi_{sc}^{0,0}}}
\newcommand\Psiscis{\operatorname{\Psi_{sc}^{0,\infty}}}
\newcommand\Psiscsi{\operatorname{\Psi_{sc}^{-\infty,0}}}
\newcommand\Psiscs{\operatorname{\Psi_{sc}^{-\infty,\infty}}}
\newcommand\Psiscalg{\operatorname{\Psi_{sc}^{\infty,-\infty}}}
\newcommand\nullHam{{\mathcal N}}
\newcommand\charD{\Sigma_{\Delta-\lambda^2}}
\newcommand\charLap{\Sigma_{\Delta-\lambda}}
\newcommand\Snl{\Sn_{\lambda}}
\newcommand\SNl{\SN_{\lambda}}
\newcommand\gammat{\tilde\gamma}
\newcommand\gammasc{\gamma}
\newcommand\Tau{\mathcal{T}}
\newcommand\taut{\tilde\tau}
\newcommand\taub{\bar\tau}
\newcommand\Nout{N^+_{\lambda}}
\newcommand\Nin{N^-_{\lambda}}
\newcommand\Nio{N^{\pm}_{\lambda}}
\newcommand\El{E_{\lambda}}
\newcommand\Elt{\tilde E_{\lambda}}
\newcommand\Eil{E^i_{\lambda}}
\newcommand\Ejl{E^j_{\lambda}}
\newcommand\Eajl{E^{\alpha_j}_{\lambda}}
\newcommand\Eilt{\tilde E^i_{\lambda}}
\newcommand\Np{N^+}
\newcommand\Nm{N^-}
\newcommand\Npm{N^{\pm}}
\newcommand\Fin{F^-(\lambda)}
\newcommand\Fini{F^-_i(\lambda)}
\newcommand\Fout{F^+(\lambda)}
\newcommand\Fouti{F^+_i(\lambda)}
\newcommand\Foutj{F^+_j(\lambda)}
\newcommand\Rout{R^+_{\lambda}}
\newcommand\Routl{R^+_{\lambda^2}}
\newcommand\Routsgnl{R^{\sign\lambda}_{\lambda^2}}
\newcommand\Rin{R^-_{\lambda}}
\newcommand\Rinl{R^-_{\lambda^2}}
\newcommand\Rinsgnl{R^{-\sign\lambda}_{\lambda^2}}
\newcommand\Rio{R^{\pm}_{\lambda}}
\newcommand\Riol{R^{\pm}_{\lambda^2}}
\newcommand\Roi{R^{\mp}_{\lambda}}
\newcommand\Roil{R^{\mp}_{\lambda^2}}
\newcommand\Riob{R^{\pm}}
\newcommand\Roib{R^{\mp}}
\newcommand\Tio{T^{\pm}}
\newcommand\Tiob{T^{\pm}_{\ff}}
\newcommand\Toi{T^{\mp}}
\newcommand\Toib{T^{\mp}_{\ff}}
\newcommand\TIiob{T_I^{\pm}}
\newcommand\Rinb{R^-}
\newcommand\Rinbsgnl{R^{-\sign\lambda}}
\newcommand\Tin{T^-}
\newcommand\Tinb{T^-_{\ff}}
\newcommand\TIinb{T^-_I}
\newcommand\Routb{R^+}
\newcommand\Routbsgnl{R^{\sign\lambda}}
\newcommand\Tout{T^+}
\newcommand\Toutb{T^+_{\ff}}
\newcommand\TIoutb{T^+_I}
\newcommand\Rlkf{(|\xib|^2-(\lambda-i0)^2)^{-1}}
\newcommand\Rlk{\rho_0(\lambda)}
\newcommand\Rmlk{\rho_0(-\lambda)}
\newcommand\Rpmlk{\rho_0(\pm\lambda)}
\newcommand\Rlka{\rho_1(\lambda)}
\newcommand\Rlkb{\rho_2(\lambda)}
\newcommand\Rilk{\rho_i(\lambda)}
\newcommand\reduced{\natural}
\newcommand\Rlf{R_0(\lambda)}
\newcommand\Rla{R_1(\lambda)}
\newcommand\Rlb{R_2(\lambda)}
\newcommand\Ril{R_i(\lambda)}
\newcommand\Rlj{R_j(\lambda)}
\newcommand\Rlft{R_0(\lambda)}
\newcommand\Rflambda{R_0^{\reduced}(\sigma)}
\newcommand\RV{R^{\reduced}_V}
\newcommand\Rfsigma{R_0^{\reduced}(\sigma)}
\newcommand\Rfsigmah{R_0^{\reduced}(\sigma^{1/2})}
\newcommand\Rfzero{R_0^{\reduced}(0)}
\newcommand\RlV{R^{\reduced}_V(\sigma)}
\newcommand\RlVi{R^{\reduced}_{V_i}(\sigma)}
\newcommand\RlVt{R_V(\lambda)}
\newcommand\RlVtL{{R}_V^L(\lambda)}
\newcommand\RlVtR{{R}_V^R(\lambda)}
\newcommand\RlVit{{R}_{V_i}(\lambda)}
\newcommand\RlVta{{R}_V^{(1)}(\lambda)}
\newcommand\RlVtk{{R}_V^{(k)}(\lambda)}
\newcommand\RlVatV{{R}_{V_{\alpha}}(\lambda)V_{\alpha}}
\newcommand\RlVatVa{{R}_{V_{\alpha_1}}(\lambda)V_{\alpha_1}}
\newcommand\RlVatVb{{R}_{V_{\alpha_2}}(\lambda)V_{\alpha_2}}
\newcommand\RlVatVk{{R}_{V_{\alpha_k}}(\lambda)V_{\alpha_k}}
\newcommand\RlVatVkk{{R}_{V_{\alpha_{k+1}}}(\lambda)V_{\alpha_{k+1}}}
\newcommand\RlVaptV{{R}_{V_{\alpha'}}(\lambda)V_{\alpha'}}
\newcommand\RlVapptV{{R}_{V_{\alpha''}}(\lambda)V_{\alpha''}}
\newcommand\RlVajtV{{R}_{V_{\alpha_j}}(\lambda)V_{\alpha_j}}
\newcommand\RlVaktV{{R}_{V_{\alpha_k}}(\lambda)V_{\alpha_k}}
\newcommand\RlVakktV{{R}_{V_{\alpha_{k+1}}}(\lambda)V_{\alpha_{k+1}}}
\newcommand\Tl{T(\lambda)}
\newcommand\Tlt{\tilde\Tl}
\newcommand\Tltp{\tilde T'(\lambda)}
\newcommand\Tltpp{\tilde T''(\lambda)}
\newcommand\Tli{T_i(\lambda)}
\newcommand\Tlit{\tilde\Tli}
\newcommand\Tlip{T_i'(\lambda)}
\newcommand\Tlipp{T_i''(\lambda)}
\newcommand\Tlj{T_j(\lambda)}
\newcommand\Tla{T_{\alpha}(\lambda)}
\newcommand\Tlaa{T_{\alpha_1}(\lambda)}
\newcommand\Tlab{T_{\alpha_2}(\lambda)}
\newcommand\Tlak{T_{\alpha_k}(\lambda)}
\newcommand\Tlakt{\tilde\Tlak}
\newcommand\Tlaj{T_{\alpha_j}(\lambda)}
\newcommand\Tlajj{T_{\alpha_{j+1}}(\lambda)}
\newcommand\Tlajp{T_{\alpha_j}'(\lambda)}
\newcommand\Tlajpt{\tilde\Tlajp}
\newcommand\Tlajt{\tilde\Tlaj}
\newcommand\Tlakk{T_{\alpha_{k+1}}(\lambda)}
\newcommand\Tlakkp{T_{\alpha_{k+1}}'(\lambda)}
\newcommand\Tlap{T_{\alpha'}(\lambda)}
\newcommand\Tlapt{\tilde\Tlap}
\newcommand\Tlapp{T_{\alpha''}(\lambda)}
\newcommand\Tkl{T^{(k)}(\lambda)}
\newcommand\Tcl{T^{\flat}(\lambda)}
\newcommand\Fl{F(\lambda)}
\newcommand\BlVt{\tilde B_V(\lambda)}
\newcommand\KBlVt{K_{\BlVt}}
\newcommand\BlVaat{B_{V_{\alpha_1}}(\lambda)}
\newcommand\BV{B_V}
\newcommand\Bone{B_1}
\newcommand\Btwo{B_2}
\newcommand\Bthree{B_3}
\newcommand\Banyj{B_j}
\newcommand\PlV{P_V(\lambda)}
\newcommand\PlVc{P_V^{\flat}(\lambda)}
\newcommand\Pl{P_0(\lambda)}
\newcommand\SVl{S_V(\lambda)}
\newcommand\Sjr{S_j^{\reduced}}
\newcommand\Rkp{{\mathcal R}^k_+}
\newcommand\Rkm{{\mathcal R}^k_-}
\newcommand\Rkpm{{\mathcal R}^k_{\pm}}
\newcommand\Phys{{\mathcal P}}
\newcommand\Pc{\overline{\mathcal P}}
\newcommand\pip{\pi^{\perp}}
\newcommand\pipa{\pi_\partial}
\newcommand\gammapa{\gamma_\partial}
\newcommand\pipah{\hat\pi_\partial}
\newcommand\pit{\tilde\pi}
\newcommand\xit{\tilde\xi}
\newcommand\zetat{\tilde\zeta}
\newcommand\etat{\tilde\eta}
\newcommand\sigmat{\tilde\sigma}
\newcommand\sigmahat{\hat\sigma}
\newcommand\thetat{\tilde\theta}
\newcommand\psit{\tilde\psi}
\newcommand\phit{\tilde\phi}
\newcommand\chit{\tilde\chi}
\newcommand\rhot{\tilde\rho}
\newcommand\xib{\bar\xi}
\newcommand\zetab{\bar\zeta}
\newcommand\thetab{\bar\theta}
\newcommand\etab{\bar\eta}
\newcommand\iotal{\iota_{\lambda}}
\newcommand\rhoat{\rhot_{\alpha_1}}
\newcommand\Lambdat{\tilde\Lambda}
\newcommand\Lambdati{\tilde\Lambda^{\text{in}}}
\newcommand\Lambdato{\tilde\Lambda^{\text{out}}}
\newcommand\Lambdatp{\tilde\Lambda^{\text{prop}}}
\newcommand\Lambdai{\Lambda^{\text{in}}}
\newcommand\Lambdao{\Lambda^{\text{out}}}
\newcommand\poles{\Lambda'}
\newcommand\rpoles{\Lambda_p}
\newcommand\thresholds{\Lambda}
\newcommand\Vt{\tilde V}
\newcommand\It{\tilde I}
\newcommand\half{{\frac{1}{2}}}
\newcommand\sigmah{\sigma^{1/2}}
\newcommand\bX{\partial X}
\newcommand\bXb{\partial \Xb}
\newcommand\Deltabt{\tilde\Delta_0}
\newcommand\strip{\Omega_T}
\newcommand\Kf{K^{\flat}}
\newcommand\Gs{G^{\sharp}}
\newcommand\Gt{\tilde G}
\newcommand\Osc{\sci\Omega}
\newcommand\OSc{{}^\Scl\Omega}
\newcommand\Osch{\sci\Omega^{\half}}
\newcommand\Oscmh{\sci\Omega^{-\half}}
\newcommand\Isc{I_{sc}}
\newcommand\os{{\text{os}}}
\newcommand\Qzl{Q^0_{-\lambda}}
\newcommand\Lie{{\mathcal L}}
\newcommand\bl{{\text b}}
\newcommand\scl{{\text{sc}}}
\newcommand\sccl{{\text{scc}}}
\newcommand\Scl{{\text{Sc}}}
\newcommand\ScLl{{\text{Sc,L}}}
\newcommand\ScRl{{\text{Sc,R}}}
\newcommand\Sccl{{\text{Scc}}}
\newcommand\sus{{\text{sus}}}
\newcommand\ssl{{\text{ee}}}
\newcommand\bzl{{\text{b0}}}
\newcommand\XXb{X^2_\bl}
\newcommand\XXbt{\Xt^2_\bl}
\newcommand\XXsc{X^2_\scl}
\newcommand\XXsct{\Xt^2_\scl}
\newcommand\XXSc{X^2_\Scl}
\newcommand\XXSct{\Xt^2_\Scl}
\newcommand\XXScL{X^2_\ScLl}
\newcommand\XXScR{X^2_\ScRl}
\newcommand\MMsc{M^2_\scl}
\newcommand\Deltab{\Delta_\bl}
\newcommand\Deltasc{\Delta_\scl}
\newcommand\DeltaSc{\Delta_\Scl}
\newcommand\DeltaScL{\Delta_\ScLl}
\newcommand\DeltaScR{\Delta_\ScRl}
\newcommand\prs{\sigma}
\newcommand\Nsc{N_\scl}
\newcommand\Nscp{N_{\scl,p}}
\newcommand\Nff{N_{\ff}}
\newcommand\Nffz{N_{\ff,0}}
\newcommand\Nffzp{N_{\ff,0,p}}
\newcommand\Nffl{N_{\ff,l}}
\newcommand\Nffml{N_{\ff,-l}}
\newcommand\Nmf{N_{\mf}}
\newcommand\Nmfz{N_{\mf,0}}
\newcommand\Nmfl{N_{\mf,l}}
\newcommand\Nmfml{N_{\mf,-l}}
\newcommand\ffb{\operatorname{bf}}
\newcommand\Ffb{\operatorname{bf'}}
\newcommand\ffsc{\operatorname{sf}}
\newcommand\ffSc{\operatorname{sf_C}}
\newcommand\Ffsc{\operatorname{sf'}}
\newcommand\rff{\rho_{\ff}}
\newcommand\rmf{\rho_{\mf}}
\newcommand\rffb{\rho_{\ffb}}
\newcommand\rffsc{\rho_{\ffsc}}
\newcommand\rFfsc{\rho_{\Ffsc}}
\newcommand\rffSc{\rho_{\ffSc}}
\newcommand\rinf{\rho_{\infty}}
\newcommand\CL{C_L}
\newcommand\CR{C_R}
\newcommand\betab{\beta_\bl}
\newcommand\betasc{\beta_\scl}
\newcommand\betaSc{\beta_\Scl}
\newcommand\BetaSc{\bar\beta_\Scl}
\newcommand\betaScL{\beta_\ScLl}
\newcommand\betaScR{\beta_\ScRl}
\newcommand\ScT{{}^\Scl T^*}
\newcommand\SccT{{}^\Scl \bar T^*}
\newcommand\ScS{{}^\Scl S^*}
\newcommand\Tb{{}^\bl T}
\newcommand\Tss{{}^\ssl T}
\newcommand\Tsc{{}^\scl T}
\newcommand\TSc{{}^\Scl T}
\newcommand\CSc{C_\Scl}
\newcommand\Lambdasc{{}^\scl\Lambda}
\newcommand\XXXb{X^3_\bl}
\newcommand\XXXsc{X^3_\scl}
\newcommand\XXXSc{X^3_\Scl}
\newcommand\XXXScO{X^3_{\Scl,O}}
\newcommand\XXXScF{X^3_{\Scl,F}}
\newcommand\XXXScS{X^3_{\Scl,S}}
\newcommand\XXXScC{X^3_{\Scl,C}}
\newcommand\KDsc{\operatorname{KD^{\half}_\scl}}
\newcommand\KDSc{\operatorname{KD^{\half}_\Scl}}
\newcommand\KDScEF{\operatorname{KD^{E,F}_\Scl}}
\newcommand\Oh{\operatorname{\Omega^{\half}}}
\newcommand\WFSc{\WF_\Scl}
\newcommand\WFtSc{\WF_{\text 3sc}}
\newcommand\WFScmf{\WF_{\Scl,\mf}}
\newcommand\WFScff{\WF_{\Scl,\ff}}
\newcommand\WFScs{\WF_{\Scl,\prs}}
\newcommand\WFScp{\WF'_\Scl}
\newcommand\WFScmfp{\WF'_{\Scl,\mf}}
\newcommand\WFScffp{\WF'_{\Scl,\ff}}
\newcommand\WFScsp{\WF'_{\Scl,\prs}}
\newcommand\Diffscc{\Diff_\sccl}
\newcommand\DiffSc{\Diff_\Scl}
\newcommand\Diffss{\Diff_\ssl}
\newcommand\DiffScc{\Diff_\Sccl}
\newcommand\DiffscI{\Diff_{\scl,\text{I}}}
\newcommand\VscI{\Vf_{\scl,\text{I}}}
\newcommand\DiffsV{\operatorname{Diff}_{\sus(V)}}
\newcommand\DiffsVsc{\operatorname{Diff}_{\sus(V),\scl}}
\newcommand\DiffsVCsc{\operatorname{Diff}_{\sus(V)-C,\scl}}   
\newcommand\Psisc{\Psop_\scl}
\newcommand\Psiscc{\Psop_\sccl}
\newcommand\Psiss{\Psop_\ssl}
\newcommand\Psisch{\Psop_{\scl,h}}
\newcommand\Psiscch{\Psop_{\sccl,h}}
\newcommand\PsiSc{\Psop_\Scl}
\newcommand\PsiScph{\Psop_{\Scl,\phi}}
\newcommand\PsiScra{\Psop_{\Scl,\rho^\sharp_a}}
\newcommand\PsiScc{\Psop_\Sccl}
\newcommand\PsiSccml{\Psop^{m,l}_\Sccl}
\newcommand\PsiScxx{\Psop^{*,*}_\Scl}
\newcommand\PsiScml{\Psop^{m,l}_\Scl}
\newcommand\PsiScmz{\Psop^{m,0}_\Scl}
\newcommand\PsiScmmz{\Psop^{-m,0}_\Scl}
\newcommand\PsiSckz{\Psop^{k,0}_\Scl}
\newcommand\PsiScmmml{\Psop^{-m,-l}_\Scl}
\newcommand\Psiscmkk{\Psop^{-k,k}_\scl}
\newcommand\Psiscmmmkk{\Psop^{-m-k,k}_\scl}
\newcommand\Psiscmoo{\Psop^{-1,1}_\scl}
\newcommand\Psiscmz{\Psop^{m,0}_\scl}
\newcommand\Psiscmmz{\Psop^{-m,0}_\scl}
\newcommand\PsiSckmkl{\Psop^{km,kl}_\Scl}
\newcommand\PsiScmplp{\Psop^{m',l'}_\Scl}
\newcommand\PsiScmmpllp{\Psop^{m+m',l+l'}_\Scl}
\newcommand\Psiscml{\Psop^{m,l}_\scl}
\newcommand\PsiScid{\Psop^{0,0}_\Scl}
\newcommand\PsiSczo{\Psop^{0,1}_\Scl}
\newcommand\PsiScmii{\Psop^{-\infty,\infty}_\Scl}
\newcommand\PsiScmiz{\Psop^{-\infty,0}_\Scl}
\newcommand\PsiScmoo{\Psop^{-1,1}_\Scl}
\newcommand\PsisCid{\Psop^{0,0}_{\scl-C}}
\newcommand\PsisC{\Psop_{\scl-C}}
\newcommand\Psiinf{\Psop_{\infty}}
\newcommand\Psiinfid{\Psop_{\infty}^0}
\newcommand\PsiFinf{\Psop_{\infty-\Fr}}
\newcommand\PsisVscml{\Psop^{m,l}_{\sus(V),\scl}}
\newcommand\PsisVsc{\Psop_{\sus(V),\scl}}
\newcommand\PsisVpsc{\Psop_{\sus(V_p),\scl}}
\newcommand\PsisVCSc{\Psop_{\sus(V)-C,\scl}}
\newcommand\SFinf{S_{\infty-\Fr}}
\newcommand\YsVC{Y^2_{\sus(V)-C,\scl}}
\newcommand\ffYsc{\ffsc_{\sus(V)}}
\newcommand\SXC{S(X;C)}
\newcommand\Ios{I_{\text{os}}}
\newcommand\pbL{\pi^2_{\bl,{\text L}}}
\newcommand\pbR{\pi^2_{\bl,{\text R}}}
\newcommand\pscL{\pi^2_{\scl,{\text L}}}
\newcommand\pscR{\pi^2_{\scl,{\text R}}}
\newcommand\PbO{\pi^3_{\bl,{\text O}}}
\newcommand\PscO{\pi^3_{\scl,{\text O}}}
\newcommand\PScO{\pi^3_{\Scl,{\text O}}}
\newcommand\PScF{\pi^3_{\Scl,{\text F}}}
\newcommand\PScC{\pi^3_{\Scl,{\text C}}}
\newcommand\PScS{\pi^3_{\Scl,{\text S}}}
\newcommand\pScL{\pi^2_{\Scl,{\text L}}}
\newcommand\pScR{\pi^2_{\Scl,{\text R}}}
\newcommand\CLF{\CL^F}
\newcommand\CLO{\CL^O}
\newcommand\CLS{\CL^S}
\newcommand\CLC{\CL^C}
\newcommand\DeltaYb{\Delta_{\bl,Y}}
\newcommand\DeltaYsc{\Delta_{\sus-\scl}}
\newcommand\diag{\operatorname{diag}}
\newcommand\Diag{\operatorname{Diag}}
\newcommand\Vf{{\mathcal V}}
\newcommand\Vb{{\mathcal V}_{\bl}}
\newcommand\Vsc{{\mathcal V}_{\scl}}
\newcommand\Vss{{\mathcal V}_{\ssl}}
\newcommand\VSc{{\mathcal V}_{\Scl}}
\newcommand\VfI{\Vf_{\text{I}}}
\newcommand\VfIq{\Vf_{\text{I},q}}
\newcommand\scH{{}^\scl H}
\newcommand\scHg{\scH_g}
\newcommand\Hss{H_\ssl}
\newcommand\Hbz{H_\bzl}
\newcommand\xh{\hat x}
\newcommand\Yh{\hat Y}
\newcommand\Zh{\hat Z}
\newcommand\Yb{\bar Y}
\newcommand\hb{\bar h}
\newcommand\xih{\hat\xi}
\newcommand\etah{\hat\eta}
\newcommand\muh{\hat\mu}
\newcommand\mub{\bar\mu}
\newcommand\nub{\bar\nu}
\newcommand\mubh{\widehat{\bar\mu}}
\newcommand\yb{\bar y}
\newcommand\zb{\bar z}
\newcommand\ub{\bar u}
\newcommand\Qb{\bar Q}
\newcommand\Wbp{{\bar W}^\perp}
\newcommand\Wp{W^\perp}
\newcommand\Kt{\tilde K}
\newcommand\Wt{\tilde W}
\newcommand\Ut{\tilde U}
\newcommand\yt{\tilde y}
\newcommand\ut{\tilde u}
\newcommand\vt{\tilde v}
\newcommand\ft{\tilde f}
\newcommand\htil{\tilde h}
\newcommand\St{\tilde S}
\newcommand\Pt{\tilde P}
\newcommand\Rt{\tilde R}
\newcommand\qt{\tilde q}
\newcommand\Qt{\tilde Q}
\newcommand\Xb{\overline{X}}
\newcommand\lambdat{\tilde\lambda}
\newcommand\betat{\tilde\beta}
\newcommand\epst{\tilde\epsilon}
\newcommand\ep{\epsilon}
\newcommand\bt{\tilde b}
\newcommand\Xt{\widetilde X}
\newcommand\Mt{\widetilde M}
\newcommand\At{\tilde A}
\newcommand\Et{\tilde E}
\newcommand\Ht{\tilde H}
\newcommand\at{\tilde a}
\newcommand\Ct{\tilde C}
\newcommand\pih{\hat\pi}
\newcommand\Rh{\hat R}
\newcommand\Ah{\hat A}
\newcommand\Bh{\hat B}
\newcommand\Ch{\hat C}
\newcommand\Gh{\hat G}
\newcommand\Hh{\hat H}
\newcommand\Qh{\hat Q}
\newcommand\Ph{\hat P}
\newcommand\Nh{\hat N}
\newcommand\Sh{\hat S}
\newcommand\Gcal{{\mathcal G}}
\newcommand\GcalC{{\mathcal G}_C}
\newcommand\Jcal{{\mathcal J}}
\newcommand\JcalC{{\mathcal J}_C}
\setcounter{secnumdepth}{3}
\newtheorem{lemma}{Lemma}[section]
\newtheorem{prop}[lemma]{Proposition}
\newtheorem{thm}[lemma]{Theorem}
\newtheorem{cor}[lemma]{Corollary}
\newtheorem{result}[lemma]{Result}
\newtheorem*{thm*}{Theorem}
\newtheorem*{prop*}{Proposition}
\newtheorem*{conj*}{Conjecture}
\numberwithin{equation}{section}
\theoremstyle{remark}
\newtheorem{rem}[lemma]{Remark}
\theoremstyle{definition}
\newtheorem{Def}[lemma]{Definition}
\newtheorem*{Def*}{Definition}
\def\signature#1#2{\par\noindent#1\dotfill\null\\*
{\raggedleft #2\par}}
\newcommand\ess{\mbox{ess}}

\renewcommand{\theenumi}{\roman{enumi}}
\renewcommand{\labelenumi}{(\theenumi)}

\title[Analytic continuation of the resolvent] 
{Analytic continuation of the resolvent of the Laplacian
on $\SL(3)/\SO(3)$}
\author[Rafe Mazzeo and Andras Vasy]{Rafe Mazzeo and Andr\'as Vasy}
\address{R.\ M.: Department of Mathematics, Stanford University, Stanford,
CA 94305}
\email{mazzeo@math.stanford.edu}
\address{A.\ V.: Department of Mathematics, Massachusetts Institute of
Technology, MA 02139}
\email{andras@math.mit.edu}
\date{August 14, 2002}

\begin{abstract}
In this paper we continue our program of extending the methods of 
geometric scattering theory to encompass the analysis of the Laplacian
on symmetric spaces of rank greater than one and their geometric
perturbations. In our previous work \cite{Mazzeo-Vasy:Scattering}
we described the resolvent, and specifically the asymptotic behavior
of the Green's function, on $\SL(3)/\SO(3)$ using methods from three-particle
scattering. Here we extend the technique
of complex scaling to symmetric spaces 
to show that the resolvent continues analytically across the spectrum.
\end{abstract}

\maketitle

\section{Introduction}
Harmonic analysis on (both global and local) symmetric spaces is an 
enduring theme in geometric and harmonic analysis, Lie theory and parts 
of number theory. The problems arising in spectral synthesis on this 
class of spaces have been approached from many points of view. This paper 
is another step in our program to study these questions using techniques 
from modern geometric scattering theory, particularly as applied to the 
analysis of quantum $N$-body Hamiltonians. In recent work 
\cite{Mazzeo-Vasy:Scattering} we have developed machinery to give
a microlocal treatment of a parametrix construction for the resolvent
of the Laplacian on the symmetric space $\SL(3,\RR)/\SO(3,\RR)$. This
produces a complete description of the singularities, particularly
at infinity, of this resolvent; however, we treated only the off-spectrum 
behaviour, for reasons of simplicity and space. Quite recently we became
aware of the fact that a full treatment of the analytic continuation of this 
resolvent past the spectrum does not seem to appear in the literature. 
(Of course, the analytic continuation of related objects, such as
the spherical functions, has been known for some time, cf.\ 
\cite{Gangolli-Varadarajan}. We also mention the paper \cite{Wallach:Powers},
where some aspects of the continuation problem are discussed from the
point of view of explicit formul\ae\ and special functions.)
We prove 
here that the resolvent of the Laplacian on $\SL(3,\RR)/\SO(3,\RR)$
does indeed have an analytic continuation. To do this we use a
combination of the microlocal techniques from \cite{Mazzeo-Vasy:Scattering}
and the method of complex scaling, which is well-known in scattering 
theory, and for which this symmetric space setting is well adapted. 

We identify $M=\SL(3,\RR)/\SO(3,\RR)$ with the set of real positive 
definite $3$-by-$3$ matrices of determinant $1$; this is a five-dimensional 
real analytic manifold. The Killing form provides a Riemannian
metric $g$. The associated Laplacian $\Delta = \Delta_g$ gives a self-adjoint 
unbounded operator on $L^2(M,dg)$, with spectrum $[\ev_0,+\infty)$, 
$\ev_0=\frac{1}{3}$. Let $R(\ev)=(\Delta-\ev)^{-1}$ be the resolvent of 
$\Delta_g$, $\ev\nin [\ev_0,+\infty)$. Fix a point $o \in M$, which
we may as well assume is the image of the identity matrix $I$ in
the identification above. The stabilizer subgroup $K_o$ (in the
natural $\SL(3,\RR)$ action on $M$) is, of course, isomorphic to $\SO(3,\RR)$.
The Green function $G_o(\ev)$ with pole at $o$ and at eigenvalue $\ev$ is, by 
definition $R(\ev)\delta_o$. It is standard that $G_o$ lies in
the space of $K_o$-invariant distributions on $M$.  

We fix the branch of the square root function $\sqrt{}$ on 
$\Cx\setminus [0,+\infty)$ which has negative imaginary part when 
$w\in\Cx\setminus[0,+\infty)$. Let $S$ denote that part of the Riemann 
surface for $\ev \mapsto \sqrt{\ev-\ev_0}$ where we continue from 
$\ev - \ev_0 \nin [0,+\infty)$ and allow $\arg(\ev-\ev_0)$ to change by 
any amount less than $\pi$. In other words, starting in the region 
$\im\sqrt{\ev-\ev_0}<0$, we continue across either of the rays where
$\im\sqrt{\ev-\ev_0}=0$ and $\re\sqrt{\ev-\ev_0}>0$, respectively $< 0$, 
allowing the argument of $\sqrt{\ev-\ev_0}$ to change by any amount less
than $\pi/2$ (so that only the positive imaginary axis is not reached). 

\noindent We now state our main result:
\begin{thm}\label{thm:main}
With all notation as above, the Green function $G_o(\ev)$ continues 
meromorphically to $S$ as a distribution. Similarly, as an operator 
between appropriate spaces of $K_o$-invariant functions, the resolvent 
$R(\ev)$ itself has a meromorphic continuation in this region, with all
poles of finite rank. 
\end{thm}

The poles which arise in this continuation of $R(\ev)$ are called 
resonances in scattering theory. Our method does not establish
that these poles necessarily exist, and we leave open the fundamental question 
of whether they actually occur. This theorem also implies that the operator
\[
R(\ev): L^2_c(M) \longrightarrow H^2_{\mathrm{loc}}(M)
\]
has a meromorphic continuation in the same region, since this is true
for its Schwartz kernel; however, since we are no longer restricting to 
$K_o$-invariant function spaces, the poles (if they exist) will no longer 
be of finite rank! 

We expect that this theorem and its proof will generalize to symmetric spaces 
of arbitrary rank without any essentially new ideas, though of course many 
new difficulties (of a mostly combinatorial nature) will need to be handled. 
We expect to return to this soon in a subsequent paper. However, it has seemed worthwhile 
presenting this method in this simple rank-$2$ setting, where the techniques
of \cite{Mazzeo-Vasy:Scattering} are already in place.  

As already noted, a key ingredient in our proof of this continuation 
is the method of complex scaling. This is an extension of dilation 
analyticity, and was originally developed in the setting of $2$-body 
scattering by Aguilar-Combes \cite{Aguilar-Combes:Analytic} and generalized to
the many-body setting by Balslev-Combes \cite{Balslev-Combes:Spectral}.
We refer to \cite{Hislop-Sigal:Spectral} and \cite[Volume 4]{Reed-Simon}
for an exposition, and to the paper \cite{Sjostrand-Zworski:Complex}
of Sj\"ostrand and Zworski for a slightly different point of view.
The main new contribution here is an adaptation
of this method to the symmetric space setting. Roughly speaking, restricted
to the space of $K_o$-invariant (or equivalently, Weyl group invariant) 
functions on the flat, the Laplacian is a three-body type Hamiltonian with 
the walls of Weyl chambers corresponding to the collision planes. This suggests that
many-body methods are natural and appropriate for this problem (and more
generally, for analysis on symmetric spaces). Our point of view combines
a uniform elliptic theory on all of $M$ with a many-body type analysis on the flats.

To give the reader a rough idea how this works, consider the hyperbolic
space $\HH^2=\SL(2,\Real)/\SO(2,\Real)$, which may be identified with
the set of two-by-two positive definite matrices $A$ of determinant $1$.
In terms of geodesic normal coordinates $(r,\omega)$ about $o=I$,
the Laplacian is given by
\begin{equation*}
\Delta_{\HH^2}=D_r^2-i(\coth r)D_r+(\sinh r)^{-2}D_\omega^2.
\end{equation*}
Now consider the diffeomorphism $\Phi_\theta:A\mapsto A^w$, $w=e^\theta$,
on $\HH^2$, $\theta\in\Real$. 
This corresponds to dilation along the geodesics through $o$,
since these have the form $\gamma_A:s\mapsto A^{cs}$, $c>0$. Thus,
in geodesic normal coordinates, $\Phi_\theta:
(r,\omega)\mapsto(e^\theta r,\omega)$. $\Phi_\theta$ defines a
group of unitary operators on $L^2(\HH^2)$ via
\begin{equation*}
(U_\theta f)(A)=(\det D_A\Phi_\theta)^{1/2}(\Phi_\theta^*f)(A),
\ J=\det D_A\Phi_\theta=w\frac{\sinh wr}{\sinh r},\ w=e^\theta.
\end{equation*}
Now, for $\theta$ real, consider the scaled Laplacian 
\begin{eqnarray*}
(\Delta_{\HH^2})_\theta & = & U_\theta \Delta_{\HH^2} U_\theta^{-1} =
J^{1/2}\Phi_\theta^* \Delta_{\HH^2} \Phi_{-\theta}^* J^{-1/2} \\
& = & J^{1/2}(w^{-2}D_r^2-iw^{-1}\coth(wr)D_r
+(\sinh(wr))^{-2}D_\omega^2)J^{-1/2}
\end{eqnarray*}
This is an operator on $\HH^2$, with coefficients which extend analytically 
in the strip $|\im\theta|<\frac{\pi}{2}$. The square root is continued 
from the standard branch near $w>0$. 
(The singularity of the coefficients at $r=0$ is only an artifact 
of the polar coordinate representation.) Note
that $(\Delta_{\HH^2})_\theta$ and $(\Delta_{\HH^2})_{\theta'}$ are
unitary equivalent if $\im\theta=\im\theta'$ because of the group properties
of $U_\theta$. The scaled operator, $(\Delta_{\HH^2})_\theta$, is not elliptic on
all of $\HH^2$ when $0<|\im\theta|<\frac{\pi}{2}$ because for $r$ large
enough, $w^2\sinh(wr)^{-2}$ can lie in $\RR^-$. However, it is elliptic in some 
uniform neighbourhood of $o$ in $\HH^2$, and its radial part 
\begin{equation*}
(\Delta_{\HH^2})_{\theta,\mathrm{rad}}=J^{1/2}(w^{-2}D_r^2-iw^{-1}\coth(wr)D_r)
J^{-1/2},
\end{equation*}
which corresponds to its action on $\SO(2)$-invariant functions, is elliptic on 
the entire half-line $r>0$. Thus, a parametrix with compact remainder can be 
constructed for $(\Delta_{\HH^2})_{\theta,\mathrm{rad}}$, and this show that 
its essential spectrum lies in $\frac{1}{4}+e^{-2i\im\theta}[0,+\infty)$. Hence 
$((\Delta_{\HH^2})_{\theta,\mathrm{rad}}-\ev)^{-1}$ is meromorphic outside this set.
In fact, it is well known that there are no poles in this entire strip
(although there are an infinite number on $|\im\theta| = \pi/2$).

Combining this with some more standard technical facts, we are in
a position to apply the theory of Aguilar-Balslev-Combes to prove
that $((\Delta_{\HH^2})_{\mathrm{rad}}-\ev)^{-1}$, and hence
$(\Delta_{\HH^2}-\ev)^{-1}$, has an analytic continuation in $\ev$ 
across $(\frac{1}{4},+\infty)$. This is done by noting that 
for $\SO(2)$-invariant functions $f,g\in L^2(\HH^2)$ and $\theta\in\Real$,
\begin{equation*}
\langle f,((\Delta_{\HH^2})_{\mathrm{rad}}-\ev)^{-1}g\rangle
=\langle U_{\bar\theta}f,((\Delta_{\HH^2})_{\theta,\mathrm{rad}}-\ev)^{-1}
U_\theta g\rangle
\end{equation*}
by the unitarity of $U_\theta$. Now if $f$, $g$ lie in a smaller (dense)
class of functions such that $U_\theta f$ and $U_\theta g$ continue
analytically from $\theta\in\Real$, then the meromorphic continuation in $\ev$ of
the right hand side is obtained by first making $\theta$ complex 
with imaginary part of the appropriate sign, and then allowing $\ev$ to 
cross the continuous spectrum of $\Delta_{\HH^2}$ without encountering 
the essential spectrum of $(\Delta_{\HH^2})_{\theta,\mathrm{rad}}$.
Hence the left hand side continues meromorphically also. With some additional
care, one can even allow $g$ to be the delta distribution at $o$, yielding the 
meromorphic continuation of the Green's function.

From the definition of $\Phi_\theta$, letting $\theta$ be complex 
corresponds to deforming $r$ into the complex plane. In fact, we 
are actually deforming $\HH^2$ in a family of totally real submanifolds
in the complexification $\SL(2,\CC)/\SO(2,\CC)$. As a point of interest,
for $\theta \notin \RR$, the maximal Grauert tube around $\HH^2$ has radius 
$\pi/2$ and so (since $\im(e^\theta r)$ is unbounded as $r \to \infty$)
the deformed contour is never contained in this Grauert tube. The critical
angle $\pi/2$ is best explained as follows. Recall the decomposition of
the Lie algebra $\sll(2,\RR) = \frakk + \frakp$, where $\frakk = 
\sol(2,\RR)$ and $\frakp$ is the space of real $2$-by-$2$ symmetric 
matrices with trace $0$. The Lie algebra of the compact dual of $\HH^2$
sits inside $\frakg_\CC$ as $\frakk + i\frakp$, and the compact dual itself 
is $\exp(i\frakp) = \sphere^2$. The exponential map on $\frakg_\CC$ has
conjugate points only in this $i\frakp$ direction (and these must occur
since $\sphere^2$ is compact). The submanifolds 
$\Phi_\theta(\HH^2)$ avoid these conjugate points until $\arg \theta = \pi/2$. 
Notice that in fact $\Phi_{\pm i\pi/2}^* \Delta_{\HH^2}\Phi_{\mp i\pi/2}^* = 
\Delta_{\sphere^2}$  (expressed in normal polar coordinates around some point $o$). 
The poles of the analytic continuation of $(\Delta_{\HH^2}-\ev)^{-1}$ 
at this critical angle can be explained using this point of view.

The case of $M=\SL(3)/\SO(3)$ is similar, except that both the elliptic
analysis of $\Delta_\theta$ and the analysis of its radial part (which
is now an operator on a 2-dimensional flat) become
more interesting: the latter (which is the most crucial part) is where
three-body scattering appears in the picture.

The plan of this paper is
as follows. In \S 2 we review the geometry of the compactification of
$M$ and the structure of the Laplacian, and sketch (a modification of) the parametrix
construction from \cite{Mazzeo-Vasy:Scattering}. The next section
defines and establishes the main facts about complex scaling in our setting.
We also prove Theorem~\ref{thm:main} here, assuming the parametrix construction 
for the scaled resolvent, to which we turn in \S 4. The brief \S 5 contains 
a few comments about further directions and questions.

The authors are grateful to Gilles Carron, Lizhen Ji and Richard Melrose 
for helpful discussions and for encouragement. 
R.\ M.\ is partially supported by NSF grant \#DMS-0204730;
A.\ V.\ is partially supported by NSF grant \#DMS-0201092 
and a Fellowship from the Alfred P.\ Sloan Foundation.
A.\ V.\ also thanks Xue Ping Wang for his kind invitation, 
and the Universit\'e de Nantes for its hospitality, 
during a stay when the work on this paper started.

\section{The structure of $M$}\label{sec:structure}
In this section we review some of the material from 
\cite{Mazzeo-Vasy:Scattering} concerning the geometry of one particularly
useful compactification of $M$, and the structure of the Laplacian 
on it. We refer to that paper for more details. 

We first define the compactification $\olM$, which is a $\Cinf$ manifold 
with corners up to codimension two. It has two boundary hypersurfaces, 
$H_\sharp$ and $H^\sharp$, which are perhaps easiest to 
describe in terms of a natural system of local coordinates
derived from the matrix representation of elements in $M$. 
Write $A\in M$ as $A=O\Lambda O^t$, with $O\in\SO(3)$ and $\Lambda$ 
diagonal. If $\fraka$ is the set of diagonal matrices of trace $0$,
then $\Lambda \in \exp(\fraka)$. The ordering of the diagonal entries is 
undetermined, but in the region where no two of them are equal,
then we denote them as $0 < \lambda_1 < \lambda_2 < \lambda_3$ (but recall
also that $\lambda_1 \lambda_2 \lambda_3 = 1$). In this region the ratios
\[
\mu = \frac{\lambda_1}{\lambda_2}, \qquad \nu = \frac{\lambda_2}{\lambda_3}
\]
are independent functions, and near the submanifold $\exp(\fraka)$ in $M$ 
we can complete them to a full coordinate system by adding the above-diagonal 
entries $c_{12}$, $c_{13}$, $c_{23}$ in the skew-symmetric matrix $T = \log O$.
Then $H^\sharp = \{\mu=0\}$ and $H_\sharp = \{\nu=0\}$,
and this coordinate system gives the $\Cinf$ structure near the
corner $H_\sharp \cap H^\sharp$. On the other hand, in a neighborhood of 
the interior of $H_\sharp$, for example, we obtain the compactification
and its $\Cinf$ structure using the coordinate system
\[
\mu, \ s = \lambda_3^{-3/2},\ c_{12},\  c_{13},\  c_{23}.
\]
We refer to \cite{Mazzeo-Vasy:Scattering}, particularly for an explanation of
this odd choice of exponent in the second coordinate, which is necessary
in order to make the coordinate change between $(\mu,s)$ and $(\mu,\nu)$
(in the overlap with the corner region) smooth. 

The Weyl group $W = S_3$ acts on the flat $\exp(\fraka)$ by permuting
the diagonal entries in $\Lambda$. For simplicity we often identify the 
flat with $\fraka$, and in this representation $W$ is generated by Euclidean 
reflections. The fixed point sets of elements of $W$ partition $\fraka$ 
into the Weyl chambers, and the fixed point sets are called the Weyl chamber 
walls. We can compactify the flat $\fraka$ to a hexagon 
$\overline{\fraka}$ by adding the faces 
at $\mu=0$ and $\nu=0$; in this picture, the corners $\mu=\nu=0$
appear well away from the closures of the Weyl chamber walls, and the sides
of the hexagon are permuted by the action of $W$. The boundaries 
$H_\sharp$ and $H^\sharp$
are the orbits of these sides under $\SO(3)$, and adjacent sides of the hexagon lie in 
different boundary hypersurfaces. The boundary hypersurfaces are each equipped with a 
fibration, with fibers $\SL(2)/\SO(2)=\HH^2$ (and base space
$\RR {\mathbb P}^2$). For example, two 
interior points of $H^\sharp$ are in the same fiber if the sum of the 
eigenspaces of the two larger eigenvalues (whose ratio is, by assumption, 
bounded in this region) is the same. The fibrations of these boundary
hypersurfaces are compatible at the corner. This gives $\olM$ a boundary 
fibration structure, similar to (but more complicated than) the ones considered 
in \cite{RBMSpec, Mazzeo:Edge, Mazzeo-Melrose:Fibred}. 

We let $\rho_\sharp$ and $\rho^\sharp$ denote boundary defining functions
for $H_\sharp$ and $H^\sharp$. Fixing $o\in M$, we may assume that
$\rho_\sharp$ and $\rho^\sharp$ are $K_o$-invariant. The function 
$\rho=\rho_\sharp\rho^\sharp$ is a {\em total boundary defining function}.
It will be important later to have such a function which is real analytic 
on $M$, and so we fix the explicit total boundary defining function
$\rho = \left(\sum_{i \neq j}\frac{\lambda_i}{\lambda_j}\right)^{-1}$.
Note that this is globally well-defined since it is a symmetric
function of the $\lambda_i$. 

There is a Lie algebra of vector fields associated to this boundary
fibration structure, namely the edge-to-edge (or simply $\ssl$) Lie 
algebra $\Vss(\olM)$ consisting of all $\Cinf$ vector fields on $\olM$ 
which are tangent to the fibers of both boundary hypersurfaces, and in 
particular tangent to $\pa\olM$. Near the corner, for example, these
are spanned over $\Cinf(\olM)$ by the basis of sections
$\mu D_\mu$, $\nu D_\nu$, $\mu D_{c_{12}}$, $\nu D_{c_{23}}$
and $\mu \nu D_{c_{13}}$. These generate the algebra of 
$\ssl$ differential operators $\Diffss(\olM)$. 
Now, $\Vss(\olM)$ is the full set of smooth sections of a vector bundle,
$\Tss\olM$, over $\olM$; its dual bundle is $\Tss^*\olM$. There is a 
principal symbol map $\sigma_{\ssl,m}$ on $\Diffss^m(\olM)$ with range 
in the space of homogeneous polynomials of degree $m$ on $\Tss^*\olM$; 
this is an extension of the usual principal symbol map on $M$ to the 
compact space $\olM$. A differential operator $P\in\Diffss^m(\olM)$
is elliptic (in the $\ssl$ calculus) if $\sigma_{\ssl,m}(P)$ is invertible 
outside the zero section. In particular, the principal symbol of 
$\Delta$ is $|\zeta|_g^2$, $\zeta \in \Tss^*\olM$, and so $\Delta$ is
an elliptic element of $\Diffss^2(\olM)$. 

If $\calF$ is any function space, we denote by $\calF^{K_o}$ the
subspace of $K_o$-invariant elements. The restriction of $\Delta$
to any one of these $K_o$-invariant subspaces is given by its radial part, 
$\Delta_{\mathrm{rad}}$. This operator plays a key role in the study
of spherical functions, see \cite{Helgason:Groups}. Moreover,
the resolvent $(\Delta-\ev)^{-1}$ commutes with the $K_o$, and its restriction to 
$K_o$-invariant functions is naturally identified with $(\Delta_{\mathrm{rad}}
-\ev)^{-1} = R(\ev)_{\mathrm{rad}}$ and is both analytically simpler and more 
amenable to the scaling method. For this 
reason, we shall deal almost exclusively with the radial Laplacian
rather than the full Laplacian from now on. 

An explicit general formula for the radial Laplacian can be found in 
\cite[Chapter II, Proposition 3.9]{Helgason:Groups}; in our case it can be 
calculated directly, cf.\ \cite{Mazzeo-Vasy:Scattering}.
In the coordinates $\mu$ and $s$, which are valid near the interior of $H_\sharp$, 
\begin{equation}
\begin{split}
\Delta_{\mathrm{rad}} = &
\frac13\left( (\mu D_\mu)^2 - \left(\frac{\mu + \mu^{-1}}{\mu-\mu^{-1}}
- \frac{s^2(\mu-\mu^{-1})}{s^4 - s^2(\mu+\mu^{-1})+1}\right) i \mu D_\mu\right)\\
&+ \frac14 \left( (sD_s)^2 - \frac{2(s^4-1)}{s^4 - s^2(\mu+\mu^{-1})+1}
i sD_s\right);
\end{split}
\label{eq:radLap2}
\end{equation}
there are similar expressions near the interior of $H^\sharp$ and near the
corner $\mu=\nu=0$. The coefficients of this operator are singular
along $\mu=0$, which is a Weyl chamber wall. This singularity is
in a certain sense genuine; in fact $\Delta_{\mathrm{rad}}$ has polar 
coordinate type singularities along all of the Weyl chamber walls, and these 
are unavoidable because the radial Laplacian is really only 
an operator on the orbifold $K_o \setminus M = \fraka/S_3$. (A simpler
example of this phenomenon is that, written as an ordinary differential
operator on $\RR^+$, the radial part of the Laplacian on $\HH^2$ 
has a regular singularity at $r=0$.) The expression (\ref{eq:radLap2})
is also interesting because we can see from it, cf.\ the discussion preceding 
Proposition 3.5 below, that the analytic continuations of its coefficients 
are also singular along the imaginary $\mu$-axis away from $\mu=0$.
This is why our method does not provide a continuation of the resolvent 
into a sector with argument larger than $\pi/2$. 

The remainder of this section is devoted to a sketch of the parametrix 
construction for $(\Delta_{\mathrm{rad}} - \ev)^{-1}$ when $\ev \notin 
[\ev_0,+\infty)$. This is meant to provide the reader some frame
of reference; the many details we omit can all be found in 
\cite{Mazzeo-Vasy:Scattering}. We shall present this construction from
a slightly different, in fact simpler, point of view than the one 
in \cite{Mazzeo-Vasy:Scattering} in that we restrict attention from the 
very beginning to the radial Laplacian. While this could have been done in 
\cite{Mazzeo-Vasy:Scattering} too, it is absolutely necessary to do so 
when we apply this construction in \S 4 to the resolvent of the scaled 
radial Laplacian, for reasons we explain below. The trade-off is that we 
are forced to deal with the singularities of the coefficients of 
$\Delta_{\mathrm{rad}}$ at the Weyl chamber walls.

Any elliptic parametrix construction is ultimately based upon the invertibility 
of certain local models for the operator in question. Indeed, the familiar interior
parametrix construction for elliptic operators uses the inverses of
homogeneous constant 
coefficient operators obtained by freezing coefficients at any point and
inverting the principal part, i.e.\ the part
with the highest homogeneity, modulo an error as indicated below. A convenient
way to obtain these inverses in this simplest case is by conjugating by the Fourier 
transform and dividing by the resulting polynomial. For elliptic operators
this polynomial does not vanish outside a compact set, and `division'
is understood as
multiplication by the product of
the reciprocal of the polynomial, or indeed its principal part,
and a function identically 1 near
infinity, but vanishing on the compact set.
The machinery of pseudodifferential 
operators provides a mechanism to patch these local inverses together. The same idea 
can be used for various degenerate problems, such as the one we have here. Now one 
must find inverses for the model operators, not only at interior points but also 
at points at infinity. The calculus of $\ssl$ pseudodifferential operators 
now provides the formal mechanism for patching these inverses together.

There are two steps in this construction. The first one mimics the local 
interior parametrix construction, but carried out uniformly in a 
neighbourhood of the diagonal. This results in the `small calculus'
parametrix, $G_1(\ev)$, which has Schwartz kernel supported in
a uniform neighbourhood of the diagonal. The operator $Q_1(\ev)$ defined by 
$(\Delta_{\mathrm{rad}} - \ev)G_1(\ev) = I - Q_1(\ev)$
is smoothing, in fact it maps $L^2(M)^{K_o}$ to $H^m_\ssl(M)^{K_o}$
for every 
$m \in \RR$, but is not compact. In the second step, we find a correction term,
$G_2(\ev)$, for this parametrix such that $I - (\Delta - \ev)(G_1(\ev) + G_2(\ev)) =
Q_2(\ev)$ is not only smoothing, but maps into functions with uniform decay at 
infinity, and thus {\it is} compact. This requires that we solve away the 
`Taylor series' of $Q_1(\ev)$ at infinity. Even to make sense of this requires 
that $\fraka$ be compactified as a manifold with corners, as we have already 
done. The construction of $G_1(\ev)$, hence of $Q_1(\ev)$, is more robust
and is carried out in \cite{Mazzeo-Vasy:Scattering} for the full Laplacian,
rather than just its radial part; in fact, as pointed out there, this does
not really use the symmetric space structure of $M$, but only its $\ssl$ boundary 
fibration structure. In our context we are faced with some minor new difficulties
because of the singularities of $\Delta_{\mathrm{rad}}$ along the Weyl chamber 
walls, but we explain below how to deal with these. On the other hand,  
the construction of $G_2(\ev)$ is more delicate and uses certain aspects
of the symmetric space structure more strongly. 

The key observation in implementing this strategy is that in certain 
regions of the compactified flat $\overline{\fraka}$, $\Delta_{\mathrm{rad}}$ 
is modelled by the simpler operators $L_{\sharp,\mathrm{rad}}$ and 
$L^\sharp_{\mathrm{rad}}$. Here
\begin{equation}
L_{\sharp,\mathrm{rad}}=\frac{1}{4}(sD_s)^2+i\frac{1}{2}(sD_s)
+\frac{1}{3}(\Delta_{\HH^2})_{\mathrm{rad}},
\label{eq:mop}
\end{equation}
and $L^\sharp_{\mathrm{rad}}$ has a very similar expression; 
of course, $(\Delta_{\HH^2})_{\mathrm{rad}} = D_\rho^2 - i\coth \rho D_\rho$ 
is the radial ($\SO(2)$-invariant) part of the Laplacian on $\HH^2$.
The somewhat cumbersome notation here is because these operators are the radial parts
of operators $L_\sharp$ and $L^\sharp$, which are obtained by replacing 
$(\Delta_{\HH^2})_{\mathrm{rad}}$ by $\Delta_{\HH^2}$, and which are models
for the full Laplacian. 
Note that, for example, $L_{\sharp,\mathrm{rad}}$ is not $S_3$-invariant, but 
if we restrict it to act on functions supported in a neighbourhood $\calU$ of 
the closure of $H_\sharp$ which does not intersect the other Weyl chamber walls,
then we can regard it as acting on $S_3$ invariant functions on $\fraka$. 
Now we can state more carefully that if $\phi_\sharp \in \Cinf_0(\calU)$, then 
$\phi_\sharp L_{\sharp,\mathrm{rad}}$ is a good local approximation to 
$\Delta_{\mathrm{rad}}$ in the sense that
\begin{equation*}
\phi_\sharp(\Delta_{\mathrm{rad}}-L_{\sharp,\mathrm{rad}})\in\rho_\sharp\Diffss^1(
\overline{\fraka}).
\end{equation*} 
The point is that this error is of only first order and has coefficients
vanishing to an extra order at $H_\sharp$, which taken together means that 
it is relatively compact with respect to $\Delta_{\mathrm{rad}}$. 

We use these model operators in both steps of the construction. For the first 
step they are used to resolve the singularities of the radial Laplacian
along the Weyl chamber walls. For the second, their resolvents are used
to solve away the restrictions of the Schwartz kernel of $Q_1(\ev)$ at the 
faces $H_\sharp$ and $H^\sharp$ (and at the corner). We invoke here
the analysis of our earlier papers
\cite{Mazzeo-Vasy:Resolvents, Mazzeo-Vasy:Scattering} concerning
the precise structure of the resolvents of product operators, which is
relevant because $L_\sharp$ and $L^\sharp$ and their radial counterparts 
are of product type. The final error term $Q_2(\ev)$ is compact on
$L^2(M)^{K_o}$ precisely when $\ev$ is not in the spectrum of 
$L_{\sharp,\mathrm{rad}}$ and $L^\sharp_{\mathrm{rad}}$, or equivalently,
$\ev \notin [\ev_0,+\infty)$. 

\medskip

We now examine each of the two steps more closely.  

\medskip

As we already indicated, the construction of the first parametrix $G_1(\ev)$ 
is quite general, and could equally well 
be carried out in the `small $\ssl$ pseudodifferential calculus' for the
resolvent of the full Laplacian. However, because we are using the
radial Laplacian, we must show how to handle the singularities of this
operator at the Weyl chamber walls. We shall rely heavily on the fact that
$\Delta_{\mathrm{rad}}$ is induced by a quotient space construction.
More specifically, the space $\fraka/S_3$ is naturally identified with 
$K_o \setminus M$ (which is the double coset space $K \setminus G / K$), 
and the operator corresponding to $\Delta_{\mathrm{rad}}$ 
in this identification lifts to the nonsingular operator $\Delta$ on $M$. 
Therefore, near the fixed point set of $K_o$ we shall use local parametrices 
for $\Delta$ which are invariant by the local isotropy group, and which 
hence push down to the quotient.  

\begin{figure}[ht]
\begin{center}
\mbox{\epsfig{file=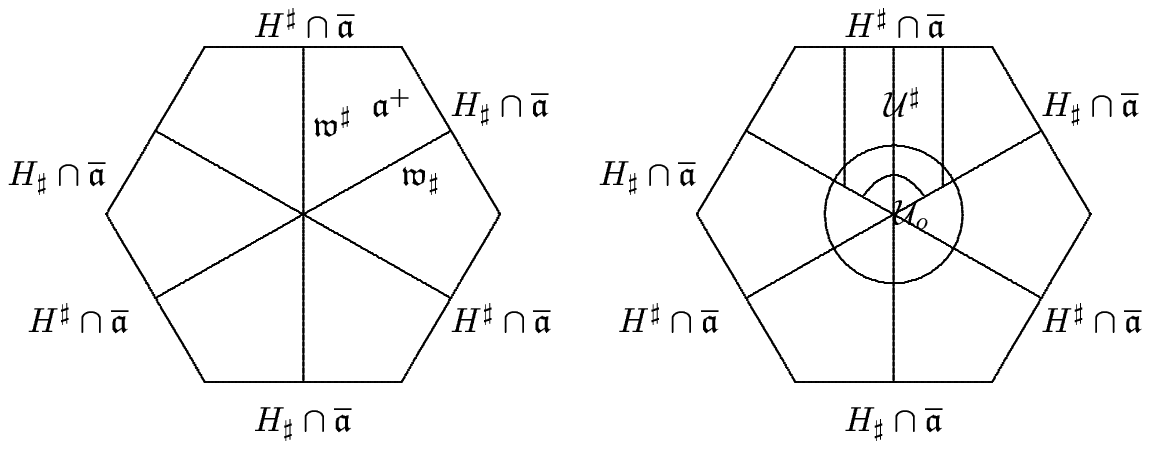}}
\end{center}
\caption{The closure of $\fraka$ in the compactifications
$\bar M$, and the neighborhoods $\calU_o$ of $0$ and $\calU^\sharp$
of the wall $\frakw^\sharp$.}
\label{fig:flat4}
\end{figure}

Let us formalize this more carefully. We divide the compactified
flat $\overline{\fraka}$ into several different regions as follows.
First choose an $S_3$-invariant ball $B_T(0) \subset \fraka$,
which we identify with a $K_o$-invariant neighbourhood around
$o$ in $M$. This corresponds to the fixed points with largest
isotropy group. Next, fix two simple roots, which we declare
to be positive, and denote $\alpha_\sharp$, $\alpha^\sharp$.
These determine two `positive' walls $\frakw_\sharp = 
\alpha_\sharp^{-1}(0) \cap (\alpha^\sharp)^{-1}(\RR^+)$ and 
$\frakw^\sharp = (\alpha^\sharp)^{-1}(0) \cap (\alpha_\sharp)^{-1}(\RR^+)$,
which together enclose the positive Weyl chamber $\fraka^+$. 
Denote by $\fraka_{+,\sharp}$ the union of $\fraka^+$, its reflection across 
the wall $\frakw_\sharp$, and the portion of this wall between them, and similarly
for $\fraka^\sharp_+$. We also let
$\overline{\fraka_{+,\sharp}}=\fraka_{+,\sharp}\cup H_\sharp$, and
similarly for $\overline{\fraka_+^\sharp}$. 
Thus $\calU^\sharp =((\alpha^\sharp)^{-1}((-T,T))
\cap \overline{\fraka_+^\sharp}) \setminus B_{T/2}(0)$ is an $S_2$-invariant neighborhood 
of $\frakw^\sharp$ which is disjoint from $0$ and all other walls, and so we 
can identify $S_3$-invariant functions on $\calU^\sharp$ with
$\SO(2)$-invariant functions on a neighborhood of
$\overline{(T/2,+\infty)}\times o'$
in $\overline{(T/2,+\infty)}\times\HH^2$
(where $o'$ corresponds to $I$ in the matrix representation of $\HH^2$).
We define the neighbourhood $\calU_\sharp$ similarly. Finally, let
$\calU_+$ be the $S_3$-orbit of the the closure of $\fraka_+$ in
$\overline{\fraka}$ minus the walls $\frakw^\sharp$ and $\frakw_\sharp$. 
This is summarized in Figure \ref{fig:flat4}. 

We lift each of these neighbourhoods to a different model space,
and in each model space we have a different pseudodifferential calculus
with which to construct parametrices. Thus $\calU_0$ lifts to a ball 
of finite radius in $M$, and we can use ordinary pseudodifferential
operators here. $\calU^\sharp$ and $\calU_\sharp$ both lift to
the model spaces $\HH^2 \times [0,1)$, where we can use the edge
calculus, and finally, we identify $\calU_+$ with itself in
$\overline{\fraka}$ and use the $\ssl$ calculus here.
This now provides the means to treat `orbifold-elliptic'
operators on the quotient $K_o\setminus M$. 

\begin{Def}\label{Def:rad-elliptic}
Suppose that $P\in\Diffss^m(M)$ is $K_o$-invariant. We say that
it is radially elliptic if its restriction to each of these neighbourhoods,
$\calU_0$, $\calU^\sharp$, $\calU_\sharp$ and $\calU_+$, 
induces operators $P_0$, $P^\sharp$, $P_\sharp$ and $P_+$, which
are symbol-elliptic in the corresponding pseudodifferential calculi
on each of the appropriate model spaces.
\end{Def}

Let us spell this out in more detail for $\Delta$. It is radial
elliptic because:

\begin{enumerate}
\item
$\Delta$ is elliptic in a neighbourhood of $o \in M$;
\item
$\frac{1}{4}(sD_s)^2+i\frac{1}{2}(sD_s) +\frac{1}{3}\Delta_{\HH^2}$ is 
elliptic in the edge calculus near $[0,1)_s\times o'$; 
\item
On $\overline{\fraka_+}$ away from the walls, 
\[
\Delta_{\mathrm{rad}} = \frac13\left((\mu D_\mu)^2+(\nu D_\nu)^2-(\mu D_\mu)(\nu D_\nu)
+i(\mu D_\mu)+i(\nu D_\nu)\right) + E,
\]
where $E$ is first order combination of the vector fields 
$\mu D_\mu$, $\nu D_\nu$ (and in addition vanishes
at the corner), cf.\ (2.5) in \cite{Mazzeo-Vasy:Scattering},
and this operator is elliptic in the $\ssl$ pseudodifferential calculus
(which in this region reduces to the $b$-calculus on a manifold with corners).
\end{enumerate}

Now, let $\{\psi_o,\psi_\sharp,\psi^\sharp, \psi_+\}$ be a smooth partition of 
unity on $\fraka$ associated to this covering, such that each function is
$S_3$-invariant.  We also select $S_3$-invariant cutoffs $\{\chi_o,\chi_\sharp,
\chi^\sharp,\chi_+\}$, such that each element is identically $1$ on a neighborhood 
of the corresponding element in the previous collection. Now suppose that
$G_o$ is a $K_o$-invariant parametrix for $\Delta-\ev$ on $B_T(0)$,
$G_\sharp = G^\sharp$ is an $\SO(2)$-invariant parametrix for $\frac{1}{4}(sD_s)^2+
i\frac{1}{2}(sD_s)+\frac{1}{3}\Delta_{\HH^2}$ on $\HH^2 \times [0,1)_s$,
and finally, $G_+$ is a parametrix $\Delta_{\mathrm{rad}}$ on a neighborhood of 
$\supp\psi_+$ in $\fraka_+$. We regard each of these operators, suitably
localized, as acting on $K_o$-invariant functions on $M$. Now set
\begin{equation*}
G_1(\ev)=\chi_o G_o(\ev)\psi_o+\chi_\sharp G_\sharp(\ev)\psi_\sharp
+\chi^\sharp G^\sharp(\ev)\psi^\sharp+\chi_+ G_+(\ev)\psi_+.
\end{equation*}
Then
\begin{equation*}
(\Delta_{\mathrm{rad}}-\ev) G_1(\ev)=\Id- Q_1(\ev)
\end{equation*}
where $Q_1(\ev):L^2(M)^{K_o}\to \Hss^m(M)^{K_o}$ for any $m$. 

We now turn to the construction of the correction term. A rather crucial
point was omitted in our earlier discussion of this second step in
the construction. Namely, while $L_{\sharp,\mathrm{rad}}$ and 
$L^\sharp_{\mathrm{rad}}$ model $\Delta_{\mathrm{rad}}$ on all of 
$\overline{H_\sharp}$, resp.\ $\overline{H^\sharp}$, there is `no room'
to patch their parametrices together on $\overline{\fraka}$ for this 
second step. (Note that this is not an issue in the preceding step because
we only used local parametrices for these operators in the interiors
of these faces, near the Weyl chamber walls, and the small $\ssl$ 
calculus is designed to handle this patching near the corner $H_\sharp\cap H^\sharp$.)
Thus, it turns out that we must perform two additional
operations: first, replace the defining functions $\rho_\sharp$ and 
$\rho^\sharp$ by $-1/\log \rho_\sharp$ and $-1/\log \rho^\sharp$, respectively,
and afterwards blow up the corners of $\overline{\fraka}$. The resulting
space is denoted $\widetilde{\fraka}$, and the closure of $\fraka_+$
inside it is denoted $\widetilde{\fraka_+}$. The corresponding operations on
$M$ produce the space $\widetilde{M}$, and as before, 
$K_o \setminus \widetilde{M}$ is identified with $\widetilde{\fraka}/S_3$.
The closures of the interiors of the faces $H_\sharp$ and $H^\sharp$
in $\widetilde{M}$ are denoted $\widetilde{H_\sharp}$ and
$\widetilde{H^\sharp}$, respectively. Instead of the neighbourhoods
$\calU^\sharp$ and $\calU_\sharp$, however, we need larger neighbourhoods:
$\tilde\calU^\sharp$ is an $S_2$-invariant neighbourhood of
$\widetilde{H^\sharp}\cap\widetilde{\fraka^+}$
which is disjoint from $\widetilde{H_\sharp}$,
and similarly for $\tilde\calU_\sharp$. There is no need to consider a
neighborhood of $0$, such as $\calU_o$, since we already have a smoothing
error term $Q_1(\ev)$, and we only need to gain decay at infinity to
obtain a compact error $Q_2(\ev)$.

\begin{figure}[ht]
\begin{center}
\mbox{\epsfig{file=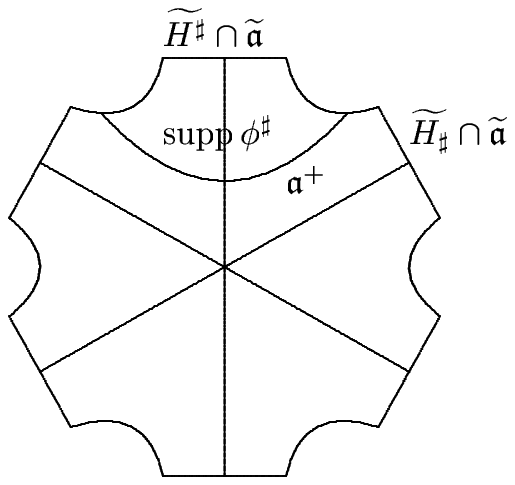}}
\end{center}
\caption{$\supp\phi^\sharp$ in $\widetilde{\fraka}$.}
\label{fig:flat3}
\end{figure}

Once these operators are lifted to this larger 
space, we can proceed much as before. Thus choose a smooth
$K_o$-invariant partitions of unity $\{\widetilde{\phi_\sharp},
\widetilde{\phi^\sharp}, \widetilde{\phi_0}\}$  and cutoffs
$\{\widetilde{\psi_\sharp}, \widetilde{\psi^\sharp}, \widetilde{\psi_0}\}$,
such that the supports of $\widetilde{\phi_\sharp}$ and 
$\widetilde{\psi_\sharp}$ are disjoint from $\widetilde{H^\sharp}$,
and similarly for $\widetilde{\phi^\sharp}$ and $\widetilde{\psi^\sharp}$,
and such that $\widetilde{\psi_\sharp}$ is identically
one on the support of $\widetilde{\phi_\sharp}$, etc. In the supports
of each of these we can use the parametrices for the respective
model operators, and so we define
\begin{equation*}
\tilde R(\ev)= \widetilde{\psi_\sharp}R_\sharp(\ev)\widetilde{\phi_\sharp}
+ \widetilde{\psi^\sharp}R^\sharp(\ev)\widetilde{\phi^\sharp},
\end{equation*}
where $R_\sharp(\ev) = (L_\sharp - \ev)^{-1}$, $R^\sharp(\ev) = 
(L^\sharp - \ev)^{-1}$.
It is not difficult to see that
\begin{equation*}
G_2(\ev)=\tilde R(\ev)Q_1(\ev)
\end{equation*}
has all the desired properties.
(We do not need to add a local parametrix for 
$\Delta-\ev$ in the support of $\psi_0$ since the error
$Q_1(\ev)$ is already smoothing.)

Namely, let
\begin{equation*}
G(\ev)=G_1(\ev)+G_2(\ev),
\ E(\ev)=(\Delta-\ev)G(\ev)-\Id,\ F(\ev)=G(\ev)(\Delta-\ev)-\Id,
\end{equation*}
and let $x\in\Cinf(M)^{K_o}$ be equal to $d(o,.)^{-1}$, the
reciprocal of the metric distance function on $M$, outside $B_1(o)$.
Then for all $m,s,l\in\RR$,
\begin{equation*}\begin{split}
G(\ev):\Hss^m(M)^{K_o}\to \Hss^{m+2}(M)^{K_o},\\
E(\ev),F(\ev):\Hss^m(M)^{K_o}\to x^s\Hss^l(M).
\end{split}\end{equation*}
The first line here follows from
\cite[Corollary~4.8]{Mazzeo-Vasy:Scattering} with $\alpha=0$, and the
second part from the proof of \cite[Proposition~5.1]{Mazzeo-Vasy:Scattering}
(which in turn uses \cite[Corollary~4.8]{Mazzeo-Vasy:Scattering}),
yielding in fact exponential decay (i.e.\ $x^s$ may be replaced by
$e^{-\beta/x}$ for sufficiently small $\beta>0$).
Thus, $E(\ev)$, $F(\ev)$ are compact operators
on $L^2(M)^{K_o}$. Hence, analytic
Fredholm theory implies that $(\Delta_{\mathrm{rad}}-\ev)^{-1}$ is meromorphic
on $\Cx\setminus[\frac{1}{3},+\infty)$
first as a bounded operator on $L^2(M)^{K_o}$, and then with the same
mapping properties as $G(\ev)$ by parametrix identities, provided it
is invertible for some value of $\ev$.
The self-adjointness and positivity of $\Delta$ in fact show {\em a priori}
that there are no
poles in $\Cx\setminus[0,+\infty)$.

Later we shall apply this construction, essentially verbatim, 
to the complex scaled version of the radial Laplacian, 
$\Delta_{\theta,\mathrm{rad}}$. The reason we have emphasized the 
use of the radial Laplacian above is that the complex scaling of the 
full Laplacian, $\Delta_{\mathrm{rad}}$, is {\it not} elliptic when 
$\theta \notin \Real$, but fortunately ellipticity is maintained for 
the radial part. We have already seen this phenomonon in our discussion 
of the scaling of $\Delta_{\HH^2}$ in the introduction. We note
also that the model operators for $\Delta_{\theta,\mathrm{rad}}$
are precisely the scaled versions $L_{\sharp,\theta,\mathrm{rad}}$
and $L^\sharp_{\theta,\mathrm{rad}}$ of the radial parts of 
$L_\sharp$, $L^\sharp$, respectively. The construction works 
just as above so long as $\ev$ is not in the spectrum of these rescaled
model operators. 

\section{Complex scaling}
We now introduce the complex scaling for $\Delta$ and $\Delta_{\mathrm{rad}}$.
This procedure produces a holomorphic family of operators for which
the essential spectrum is shifted. The parametrix construction outlined
at the end of the last section is then sufficient to produce a good
parametrix, and this leads ultimately to the analytic continuation
of the resolvent of $\Delta$. 

Let us first review the structure of complex scaling. Fix a 
connected simply connected domain $D\subset\Cx$.  We wish to
define a family of operators $U_\theta$, $\theta\in D$, 
and dense subspace $\calA\subset \Hss^\infty(\olM)^{K_o}$, 
with the following properties:
\begin{enumerate}
\item
If $\theta\in D\cap\Real$, $U_\theta$ is unitary on $L^2(M)$ and is bounded 
on all Sobolev spaces, and furthermore, $U_0=\Id$;
\item
For $f\in\calA$, $\theta \to U_\theta f$ has an analytic extension 
from $D \cap \Real$ to all of $D$ with values in $L^2(M)^{K_o}$;
\item
The subspace $U_\theta\calA$ is dense in $L^2(M)^{K_o}$.
\end{enumerate}

Now, given such a family $U_\theta$, we wish to define the family of operators 
$\Delta_\theta=U_{\theta}\Delta U_{\theta}^{-1}$. When $\theta\in\Real$,
these are defined in the obvious way since $U_\theta$ is bounded, and
they are all unitarily equivalent. The $\Delta_\theta$ are differential
operators on $M$ with coefficients that are analytic in $\theta$ in
$|\im\theta|<\frac{\pi}{2}$, as we
show below,
hence for e.g.\ $\Cinf_c(M)$, $\Delta_\theta f$ is analytic in $\theta$.
This gives an analytic continuation of $\Delta_\theta$ from $\theta\in\Real$.
The precise statement that we actually need is that $\Delta_\theta$ is
analytic of type A, see Proposition~\ref{prop:type-A}.

In complex scaling, one uses the resolvent $(\Delta_\theta-\ev)^{-1}$
of the scaled operator $\Delta_\theta$ to extend $R(\ev)=(\Delta-\ev)^{-1}$.
More precisely, we show that
the matrix elements $\langle f,R(\ev)g\rangle$ continue to all $\theta \in D$
so long as $f, g\in\calA$; using the density of $\calA$ this suffices
for purposes of spectral theory. In particular, we have
$(\Delta_\theta - \ev)^{-1} = U_\theta R(\ev) U_{\theta}^{-1}$. This
equality is valid initially when $\ev$ is in the resolvent set common
to both operators and $\theta$ is real, but after we construct the resolvent
$(\Delta_\theta-\ev)^{-1}$
by parametrix methods, this will give the meromorphic continuation of $R(\ev)$.

We define a family of unitary operators $U_{\theta}$, $\theta\in\Real$ on 
$L^2(M,dg)$ as follows. First write $w = e^\theta >0$ and let $\Phi_{\theta}:
M\to M$ be the diffeomorphism given by $A\mapsto A^w$. Then set
\begin{equation*}
(U_{\theta}f)(A)=(\det D_A\Phi_{\theta})^{1/2} f(\Phi_{\theta}(A)).
\end{equation*}
Hence $U_{\theta}$ is pull-back by $\Phi_{\theta}$ multiplied by 
the square root of the Jacobian determinant of $\Phi_{\theta}$; this
factor is included precisely to make $U_{\theta}$ unitary on $L^2(M)$. 
Clearly $U_{\theta}$ preserves $L^2(M)^{K_o}$ and is unitary on 
this subspace.

This definition is motivated as follows. For $A \in M$, $B=\log A$
is well defined (as the principal branch of the logarithm of a positive definite 
matrix). The scaling $A \mapsto A^w$ corresponds to 
the more familiar dilation $B\mapsto wB$, which plays an important
role in Euclidean scattering.

Note that (not necessarily unit speed) geodesics through $o=I$ are given by
$\gamma_A:\Real\ni s\mapsto A^{cs}$, $c>0$. Thus, for $\theta\in\Real$,
$\Phi_\theta$ is dilation along geodesics through $o$. Thus, in geodesic
normal coordinates $(r,\omega)$ around $o$, $\Phi_\theta:(r,\omega)\mapsto
(e^\theta r,\omega)$. The map $A\mapsto A^w$ is well defined for all 
complex $\theta$, and it corresponds to deforming $r$ into a sector of the 
complex domain. As for the case of $\HH^2$, we are rotating 
$M$ into $\SL(3,\CC)/\SO(3,\CC)$, but only to angles which avoid
the compact dual $\exp(\frakk + i\frakp)$, cf. the corresponding
discussion for $\HH^2$ in the introduction (also for the language used
here).  

Although $\Phi_\theta$ is a real analytic map on $M$, its regularity 
properties on $\olM$ are slightly more subtle. The following is easy to
verify in local coordinates:
\begin{lemma}
For $\theta\in\Real$,
$\Phi_\theta$ extends to a homeomorphism of $\olM$;
moreover, this extension is a `conormal diffeomorphism' 
in the sense that $\Phi^*_\theta:S^m(\olM)\mapsto S^{m w}(\olM)$, 
where $w=e^\theta$ as usual and $S^m(\olM)$ denotes the symbol space.
If $C$ denotes a matrix representative of an element of the fiber
$\HH^2=\SL(2)/\SO(2)$ on one of the boundary hypersurfaces, then 
$\Phi_\theta$ is natural in that it induces the map 
$C\mapsto C^w$.
\end{lemma}

\begin{cor}
The Jacobian determinant satisfies
\begin{equation*}
(\det D_A\Phi_{\theta})^{1/2}\in \rho^{-w+1}S^0(\olM)^{K_o},\ w=e^\theta,
\end{equation*}
and it extends analytically to the region
\begin{equation*}
\{\theta\in\Cx:\ |\im\theta|<\frac{\pi}{2}\}.
\end{equation*}
\end{cor}

\begin{rem} Another standard way of constructing the family $U_\theta$
is via a family of diffeomorphisms generated by some choice of
vector field particularly well-adapted to the operator. The 
work of Christian G\'erard \cite{Gerard:Distortion} contains a good
example of this approach. Our choice was dictated by its simplicity
and naturality given the matrix structure. In fact, by the general theory 
of complex scaling, cf.\ \cite[Theorem~16.4]{Hislop-Sigal:Spectral}, the 
final analytic continuation of the resolvent does not depend on the choice.
\end{rem}

For later purposes we define another family of unitary 
operators $U_{\theta,T}$, associated to a family of diffeomorphisms 
$\theta \to \Phi_{\theta,T}$ of $M$, where each 
$\Phi_{\theta,T}$ is the identity in the ball $B_T(o)$ and equals 
$\Phi_{\theta}$ outside a compact subset of $M$, and 
which depends analytically on $\theta$. 

To do this, fix $T>0$ and a nondecreasing cutoff function $\phi\in\Cinf(\Real;[0,1])$ 
which is identically $1$ near $\infty$ and identically $0$ on $[0,T]$. 
Now let $r(A)=d(o,A)$, where $d$ is the metric distance, and define
\begin{equation*}
\Phi_{\theta,T}(A)=A^{1+(w-1)\phi(r)},\ w=e^\theta.
\end{equation*}
Then $\Phi_{\theta,T}(A)=A$ if $d(o,A)\leq T$, and $\Phi_{\theta,T}(A)=
A^w$ outside a compact subset of $M$. Moreover, $\Phi_{\theta,T}(A)$ is
analytic in $\theta$. We must show that $\Phi_{\theta,T}$
is a diffeomorphism when $\theta$ is real and near $0$, and that it extends
analytically to complex $\theta$. 

\begin{lemma}
There exists $\delta>0$ such that $\Phi_{\theta,T}:M\to M$ is a
diffeomorphism when $\theta \in \Real$, $e^\theta>1-\delta$. In addition, 
for any $A$, $(\det D_A\Phi_{\theta,T})^{1/2}$ extends analytically to
the region 
\begin{equation*}
\{\theta\in\Cx:\ |\im\theta|<\frac{\pi}{2},\ e^\theta\nin(-\infty,1-\delta)\}.
\end{equation*}
\end{lemma}

\begin{proof}
First note that geodesics through $o$ are given by $\{A^s:\ s\in\Real\}$,
and $\Phi_{\theta,T}$ preserves these. In fact, in terms of geodesic normal 
coordinates $(r,\omega)\in [0,+\infty)\times \sphere^4$ centered at $o$, 
$\Phi_{\theta,T}$ is given by
\begin{equation*}
(r,\omega)\mapsto ((1+(w-1)\phi(r))r,\omega).
\end{equation*}
This is a diffeomorphism provided the derivative of the first component 
is nonvanishing, i.e.
\begin{equation*}
1+(w-1)(\phi(r)+r\phi'(r)) \neq 0.
\end{equation*}
But $\phi(r)+r\phi'(r)$ is smooth and equals $1$ for $r$ large,
hence is bounded; in particular, its range is contained in some 
interval $[0,\delta^{-1}]$. Thus $\Phi_{\theta,T}$ is a diffeomorphism 
provided $-\frac{1}{w-1}$ is {\em not} in this range, i.e.\ if
$w\nin (-\infty,1-\delta]$. This proves the first claim.

The second claim follows immediately from the same calculation. 
We note that the limitation $|\im\theta|<\frac{\pi}{2}$ arises
because in these coordinates the metric is given by $dr^2+h(r,\omega,d\omega)$, 
and the determinant of $h$ extends to be holomorphic and nonvanishing in this strip.
\end{proof}

As before, set
\begin{equation}
(U_{\theta,T}f)(A)=(\det D_A\Phi_{\theta,T})^{1/2} f(\Phi_{\theta,T}(A)).
\label{eq:defuo}
\end{equation}

The subspace $\calA$ is constructed below in Proposition~\ref{prop:dense}, 
and we show there that it has the desired properties for both 
$U_\theta$ and $U_{\theta,T}$.

In any case, because of the geometric nature of this specific choice of $U_\theta$,
we may define the family of differential operators $\Delta_{\theta}=U_{\theta}
\Delta U_{\theta}^{-1}$, without worrying about functional analytic issues of domain. 
These are all $K_o$-invariant on $M$ and have coefficients which depend analytically on 
$\theta$ in some strip containing the real axis. Examining the expression 
(\ref{eq:radLap2}), we see that the coefficients there extend holomorphically to 
$|\im\theta|<\frac{\pi}{2}$.
The delicate term is 
$\frac{\mu^{-1}+\mu}{\mu^{-1}-\mu} \, \mu D_\mu$. Note that the apparent singularity 
at $\mu=1$ is only a polar coordinate singularity, as can be seen by setting 
$r=-\log\mu$, which transforms this term to $-\coth r\, D_r$. However,
in the complex scaling, $r$ is replaced by $wr$, $w=e^\theta$, and when 
$w=i$ ($\theta=i\,\frac{\pi}{2}$), $\coth ir = i \cot r$, which has
(non locally integrable) singularities at $r = \pi, 2\pi, \ldots$.
This shows that $|\im\theta|=\frac{\pi}{2}$ is a genuine boundary 
for the continuation of this operator. 

We have established the smoothness of the coefficients of $\Delta_\theta$ 
when $|\im\theta|<\frac{\pi}{2}$. However, we also need 
information about its behaviour at $\pa \olM$, and we turn to this now.

When $\theta\in\Real$, the pullbacks of the boundary defining functions 
$\mu$ and $\nu$ by $\Phi_\theta$ are $\mu^w$ and $\nu^w$, and similarly
the total boundary defining function $\rho$ pulls back to $f \rho^w$
where $f$ is some smooth nonvanishing function on $\olM$. We then compute
the pullbacks of the generating vector fields for the structure algebra 
$\Vss$. The vector fields $\mu D_\mu$ and $\nu D_\nu$ are just multiplied
by a factor of $w$ in this pullback, but  $\mu D_{c_{12}}$, 
$\nu D_{c_{23}}$ and $\mu\nu D_{c_{13}}$ are transformed to
to $\mu^w D_{c_{12}}$, $\nu^w D_{c_{23}}$ and $\mu^w\nu^w D_{c_{13}}$,
respectively, and these are no longer smooth on $\olM$. This could
be remedied by changing the smooth structure, but we do not need to
consider this since we are only interested in the radial part of 
$\Delta_\theta$. To see this, suppose we are checking radial ellipticity of
the scaled operator near $\mu=1$; then $\mu^w D_{c_{12}}$ is smooth here,
while $\nu^w D_{c_{23}}$ and $\mu^w\nu^w D_{c_{13}}$ drop out. 

\begin{prop}\label{prop:type-A}
If $\theta\in\Cx$ has $|\im\theta|<\frac{\pi}{2}$, then $\Delta_\theta$
is radially elliptic. The coefficients of this operator
are conormal, but not one-step polyhomogeneous. Moreover,
$\theta \to \Delta_\theta$ is an analytic type-A family on $L^2(M)^{K_o}$
with domain $\Hss^2(\olM)^{K_o}$.
\end{prop}

\begin{proof}
The first part is easy from the explicit formula.
Now, since $\Delta_\theta$ is radially elliptic, the domain of
$\Delta_{\theta,\mathrm{rad}}$ is $\Hss^2(\olM)^{K_o}$. For any
$f\in\Hss^2(\olM)^{K_o}$, the map $\theta\mapsto\Delta_\theta f\in 
L^2(M)^{K_o}$ is strongly analytic, and this is what it means for 
$\Delta_\theta$ to be an analytic family of type A.
\end{proof}

We now state the conclusion of Proposition~\ref{prop:ess-spec} below,
that for $|\im\theta|<\pi/2$ and as an operator on $L^2(M)^{K_o}$,
\begin{equation}\label{eq:ess-spec-8}
\spec_{\mathrm{ess}}(\Delta_{\theta})=
\ev_0+e^{-2i\im\theta}[0,+\infty).
\end{equation}
In other words, the complex scaling rotates the essential spectrum of 
$\Delta$ around the threshold at the bottom of the spectrum by an angle 
$2\im\theta$. The proof of this key fact in \S 4 uses simple three-body 
techniques, and is not hard using our geometric approach; we suspect that
this is likely to be much less tractable with other, more traditional
approaches. Granting this result, we now proceed with the rest of the argument.

We shall write $R(\ev,\theta)$ for $(\Delta_{\theta,\mathrm{rad}} - \ev)^{-1}$.
To be definite, in the remainder of this section we shall treat only
the analytic continuation of $R(\ev)_{\mathrm{rad}} = R(\ev,0)$ from the lower right 
quadrant $\im(\ev-\ev_0)<0$ (and through the ray $(\ev_0,+\infty)$); the continuation 
of $R(\ev)$ from $\im(\ev-\ev_0)$ is handled nearly identically. 

Normalize so that $\arg(\ev-\ev_0)\in(-\pi,0)$ for $\im(\ev-\ev_0)<0$.
Fix any $\beta \in  (0,\frac{\pi}{2})$, and define the Riemann surface
\begin{equation*}
S_\beta=\{\ev\in S:\ -\pi<\arg(\ev-\ev_0)<2\beta\}.
\end{equation*}
Thus for $\ev \in S_\beta$, $-\frac{\pi}{2}<\arg\sqrt{\ev-\ev_0}<\beta$.

The main point, roughly speaking, is that when $-\frac{\pi}{2}<\im\theta<0$, 
$\Delta_{\theta}-\ev$ is a holomorphic family of operators on $S_\beta$
with values in the space of radially elliptic elements in $\Diffss^2(\olM)$.
Then, using \eqref{eq:ess-spec-8}, $R(\ev,\theta)$ 
is meromorphic in $\ev$ outside $\ev_0+e^{-2i\im\theta}[0,+\infty)$, 
with values in bounded operators on $L^2(M)^{K_o}$. This family has
only finite rank poles, and these are the poles of the continuation of 
$R(\ev)_{\mathrm{rad}}$ if we choose $\theta$ so that $\beta<-\im\theta<\frac{\pi}{2}$.

More precisely, we shall use the theorem of Aguilar-Balslev-Combes. This
theorem has three kinds of hypotheses: some about the scaling
$U_\theta$, some about the set of analytic vectors $\calA$ (these
two are summarized in (i)-(iii) at the beginning of \S 3),
and some about the scaled operators
$\Delta_\theta$. We have explained all of these ingredients except
$\calA$, and we postpone the discussion of $\calA$ until after the statement
of the theorem. There, we specify $\calA$ directly,
as a space of entire functions decaying at infinity, and show that the
hypotheses (i)-(iii) at the beginning of \S 3 are satisfied for
$D=\{\theta:\ |\im\theta|<\frac{\pi}{4}\}$. This
guarantees analytic continuation of $R(\ev)$ to $S_{\pi/4}$, rather
than $S_{\pi/2}$ as in Theorem~\ref{thm:main}, but the full strength of
Theorem~\ref{thm:main} is obtained subsequently by iterating this argument.

There is a general abstract theorem, due to Nelson, cf.\ 
\cite[Volume~2]{Reed-Simon}, which provides a dense set of
analytic vectors for generators of groups of unitary operators, essentially
by using the functional calculus, more specifically the spectral projection
to bounded subsets of $\Real$. Nelson's theorem would allow
us to replace $D=\{\theta:\ |\im\theta|<\frac{\pi}{4}\}$ by
$D=\{\theta:\ |\im\theta|<\frac{\pi}{2}\}$ directly in the following theorem,
and give the analytic continuation of $R(\ev)$ to $S_{\pi/2}$.
However, the explicit definition of the subspace $\calA$ which we use 
is intended to make the proof more transparent.

\begin{thm*} {\rm (\,}\cite[Theorem~16.4]{Hislop-Sigal:Spectral}\,{\rm )}
Suppose that $U_\theta$ and $\calA$ satisfy the hypotheses (i)-(iii) listed
in the beginning of \S 3, and that $\Delta_\theta$ is a type-A analytic 
family in the strip $D=\{\theta:\ |\im\theta|<\frac{\pi}{4}\}$, and
\eqref{eq:ess-spec-8} holds for $\theta\in D$. Then
\begin{enumerate}
\item
For $f,g\in\calA$, $\langle f,R(\ev)_{\mathrm{rad}}g\rangle$ has a meromorphic 
continuation to $S_{\pi/4}$.
\item
The poles of the continuation of $\langle f,R(\ev)_{\mathrm{rad}}g\rangle$ to $S_\beta$,
$\beta<\frac{\pi}{4}$, are the union of all eigenvalues of $\Delta_{\theta,\mathrm{rad}}$ 
for $\beta<-\im\theta<\frac{\pi}{4}$.
\item
The poles are independent of the choice of $U_\theta$ in the sense
that if $U'_\theta$ and $\calA'$ also satisfy (i)-(iii) and if
$\calA\cap\calA'$ is dense, then the eigenvalues of $U'_\theta
\Delta_{\mathrm{rad}} (U'_\theta)^{-1}$ are the same as those of 
$\Delta_{\theta,\mathrm{rad}}$.
\end{enumerate}
\end{thm*}

We briefly outline the proof of the first part since the idea is simple.
We must relate $R(\ev,\theta)$ and $R(\ev)_{\mathrm{rad}}$, so fix some
$\ep>0$, and suppose that 
\begin{equation*}
\theta \in \Omega_{\ep}=\{-\ep<\im\theta<\frac{\pi}{4}\}
\qquad \mbox{and}\qquad \arg(\ev-\ev_0)\in(-\pi,-\ep).
\end{equation*}
When $\theta$ is real, $U_{\theta}$ is unitary and so
\begin{equation*}
\langle f,R(\ev)_{\mathrm{rad}}g\rangle
=\langle U_{\bar\theta}f,(U_{\theta}R(\ev)_{\mathrm{rad}}U_{\theta}^{-1})
U_{\theta}g\rangle;
\end{equation*}
in addition $U_{\theta}R(\ev)_{\mathrm{rad}}U_{\theta}^{-1}=R(\ev,\theta)$, and hence
\begin{equation}\label{eq:anal-cont-16}
\langle f,R(\ev)_{\mathrm{rad}}g\rangle
=\langle U_{\bar\theta}f,R(\ev,\theta)U_{\theta}g\rangle.
\end{equation}
The left side of this equation is independent of $\theta$; on the other hand,
the right side is analytic in $\theta$ on $\Omega_{\ep}$, and 
hence is constant on this domain. This holds when $\arg(\ev-\ev_0)\in(-\pi,-\ep)$. 

The extension of $\langle f,R(\ev)_{\mathrm{rad}}g\rangle$ to $S_\beta$ is done as 
follows. Take $\theta$ with 
$-\im\theta\in(\beta,\frac{\pi}{4})$. For $\ev\in\Cx$ with $\im(\ev-\ev_0)<0$,
$\langle f,R(\ev)_{\mathrm{rad}}g\rangle$ is given by the right hand side of
\eqref{eq:anal-cont-16}. But this right side is analytic in $\ev$
away from the spectrum of $\Delta_{\theta,\mathrm{rad}}$, and meromorphic away from 
its essential spectrum, hence is meromorphic on $S_\beta$, as claimed.
This continuation is clearly independent of the choice of $\theta$ with
$-\im\theta\in(\beta,\frac{\pi}{4})$ since any such continuation is
a meromorphic function of $\ev$ that agrees with a given function on
an open set.

It remains to define the subspace $\calA$ of analytic vectors. Let $\fraka^{\Cx}$ 
denote the complexification of the Cartan subspace $\fraka$ (which is
identified with the space of diagonal matrices $\exp(\fraka)$ in $M$).
The action of the Weyl group $S_3$ extends naturally to this complexification.
We define $\calA$ to be the space of entire $S_3$-invariant functions $f$ on 
$\fraka^{\Cx}$ which decay faster than any power of $(1+|z|^2)^{-1/2}$ 
on every slab $\{z\in\fraka^\Cx:\ |\im z|\leq C|\re z|\}$, $0<C<1$:
for every $0<C<1$ and $N>0$,
\begin{equation*}
\sup_{|\im z|\leq C|\re z|} |f(z)|(1+|z|^2)^{N}<+\infty.
\end{equation*}

\begin{prop}\label{prop:dense}
When $|\im\theta\,|<\frac{\pi}{4}$, $U_{\theta}\calA$ is dense in
$L^2(M)^{K_o}$.
\end{prop}

\begin{proof}

The space $L^2(M)^{K_o}$ can be identified with an $L^2$ space
of $S_3$-invariant functions on $\fraka$ with respect to a 
singular measure. In fact, if $dg_0$ is the standard Euclidean 
density on $\fraka$ and $\rho$ is a $K_o$-invariant total boundary 
defining function on $\olM$, then there exists a smooth $S_3$-invariant 
function $a$ on $\fraka$ with $\rho^2 a\in\Cinf(\olM)^{K_o}$ such that 
$L^2(M,dg)^{K_o} = L^2(\fraka, a \, dg_0)^{S_3}$. However, there are 
effectively polar coordinate singularities which cause $a$ to vanish
at the Weyl chamber walls.

Let $c_t \, e^{-|x|^2/t}$ be the fundamental solution of the Euclidean 
heat equation, $c_t=1/\pi t$. Then if $f \in \calC^0_c(\fraka)^{S_3}$, we let
\begin{equation*}
f_t(x)=f\ast c_t e^{-|x|^2/t}
=\int f(y) c_t\, e^{-|x-y|^2/t}\,dy.
\end{equation*}
We claim that for $t > 0$, $f_t\in\calA$. Indeed, if $\theta \in \Real$, 
$f_t(e^\theta x)= \int f(y) c_t e^{-(e^\theta x-y)^2/t}$, 
and this continues analytically in $\theta$ because $\exp(-(z-y)^2)$
is entire in $z$ and decreases faster than any power of $(1+|z|^2)^{-1/2}$ 
in $|\im z|<C|\re z|$. Finally, since $S_3$ acts by reflections
and these commute with the heat kernel, each $f_t(x)$ is $S_3$ invariant.
This proves the claim.

To prove that $\calA$ is dense in $L^2(M,dg)^{K_o}$, we note first
that $e^{|x|^2}f_t$ is uniformly bounded when $t<1$, and
$\sup e^{|x|^2}|f(x)-f_t(x)|\to 0$ as $t\to 0$. Since
$a\, e^{-|x|^2}\in L^1(dg_0)$, we see that for such $f$, 
$f_t\to f$ in $L^2(\fraka, a \, dg_0)^{S_3}$. This is sufficient,
since the set of such $f$ is dense in $L^2(\fraka,a\, dg_0)^{S_3}$.

This argument requires only minor modifications to prove that
$U_\theta\calA$ is also dense in $L^2(M,dg)^{K_o}$. Namely,
for $f$ as above, the function 
\begin{equation*}
f_t(x)=\int f(y) c_t e^{2\theta} e^{-|x-e^\theta y|^2/t}\,dy \in \calA,
\end{equation*}
and we only need to show that $f_t(e^\theta x)\to f$
in $L^2(\fraka, a\, dg_0)^{S_3}$, or even more simply, that
$\sup e^{|x|^2}|f(x)-f_t(e^\theta x)|\to 0$ as $t\to 0$. However,
\begin{equation*}
f_t(e^\theta x)=\int f(y) c_t e^{2\theta}
e^{-e^{2\theta}|x-y|^2/t}\,dy,
\end{equation*}
and $\int c_t e^{2\theta} e^{-e^{2\theta}|x-y|^2/t}\,dy=1$, and so the 
standard argument finishes the proof.
\end{proof}

\begin{rem}
Note that $e^{-\gamma z^2}\to 0$ as $|z|\to\infty$ in $|\im z|\leq (1-\ep)
|\re z|$ for any $\ep>0$, but this fails when $\ep=0$. 
We have required $\arg e^\theta<\frac{\pi}{4}$ to ensure that $U_\theta 
(e^{-\gamma x^2})$, which is a multiple of $e^{-\gamma (e^\theta x)^2}$,
lies in $\calA$. 
\end{rem}

\begin{cor}\label{cor:dense}
For $|\im\theta\,|<\frac{\pi}{4}$, $U_{\theta}\calA$ is dense in $\Hss^m(M)
^{K_o}$ for all $m > 0$. 
\end{cor}

\begin{proof}
Implicit in the definition of these Sobolev spaces, but cf.\ 
\cite{Mazzeo-Vasy:Scattering} for an explanation, $(\Delta+1)^{m/2}:
\Hss^m(M)^{K_o}\to L^2(M)^{K_o}$ is an isomorphism. Thus, $f_t\to f$ as 
$t\to 0$ in $\Hss^m(M)^{K_o}$ if and only if $(\Delta+1)^{m/2}f_t\to 
(\Delta+1)^{m/2}f$ in $L^2(M)^{K_o}$. So given $f\in \Hss^m(M)^{K_o}$, 
let $g=(\Delta+1)^{m/2}f$. Since $\calA$ is dense in $L^2(M)^{K_o}$, 
there exists a family $g_t\in\calA$ with $g_t\to g$ as $t\to 0$ in 
$L^2(M)^{K_o}$. Now let $f_t=(\Delta+1)^{-m/2}g_t$ and note that 
$f_t\in\calA$. Thus, $f_t\to f$ in $\Hss^m(M)^{K_o}$ as desired.
\end{proof}

The continuation of $U_\theta g$ may well exist for functions
or distributions $g$ which do not lie in $\calA$. For example,
if $g=\delta_o$, the delta distribution at $o$, then its homogeneity
shows that for $\theta$ real, $U_{\theta}\delta_o=(\det 
D_o\Phi_{\theta})^{-1/2}\delta_o$. Clearly, then, $U_{\theta}\delta_o$
extends to be analytic in $\theta$ (e.g.\ with values in a Sobolev space),
and so the Green function, $R(\ev)\delta_o$ also extends 
via $\langle f,R(\ev)\delta_o\rangle$ for $f\in\calA$.

Note that this does not quite say that $R(\ev)\delta_o$ is a distribution, 
since that would require that the right hand side of 
\eqref{eq:anal-cont-16} be defined for any $f\in\Cinf_c(\fraka)^{S_3}$,
while for most $f$, $U_{\bar\theta} f$ does {\em not} have an analytic extension 
from the real axis. However, we can prove this using the deformed
group of unitary operators, $U_{\theta,T}$, defined in (\ref{eq:defuo}).
Recall that the associated diffeomorphisms $\Phi_{\theta,T}$ fix points in
$B_T(o)$ and equal $\Phi_{\theta}$ for $|x|$ sufficiently large.
We use precisely the same arguments as above to establish the density of 
$U_{\theta,T}\calA$. Hence by the 
uniqueness part of the Aguilar-Balslev-Combes theorem, the induced analytic 
extensions agree with one another no matter the value of $T$, and also agree
with the extension associated to $U_\theta$. But if $f\in\Cinf_c(B_T(o))^{S_3}$,
then $U_{\theta,T}f=f$ and so $U_{\theta,T}f=f$ has an analytic extension
to $\theta\in\Cx$. Arguing as before, the formula 
\begin{equation}\label{eq:anal-cont-32}
\langle f,R(\ev)_{\mathrm{rad}}\delta_o\rangle
=\langle U_{\bar\theta,T}f,R(\ev,\theta,T)U_{\theta,T}\delta_o\rangle
=\langle f,R(\ev,\theta,T)\delta_o\rangle
\end{equation}
shows that $R(\ev)_{\mathrm{rad}}\delta_o$ does indeed extend analytically 
as a distribution to $S_{\pi/2}$ since the right hand side has this
property.

Although we have only constructed a subset $\calA \subset L^2(M)^{K_o}$
for which $U_{\theta}\calA$ is dense in $L^2(M)^{K_o}$ when
$|\im\theta|<\pi/4$, we can still continue $R(\ev)$ to $S_{\pi/2}$,
for which the formula \eqref{eq:anal-cont-16} requires
larger $\im\theta$.

\begin{thm*}[Theorem~\ref{thm:main}]
The Green function $G_o(\ev)$ continues 
meromorphically to $S_{\pi/2}$ as a distribution.
\end{thm*}

\begin{proof}
We have shown that the hypotheses of the Aguilar-Balslev-Combes theorem are
satisfied for $D=\{\theta:\ |\im\theta|<\frac{\pi}{4}\}$ (for either
$U_\theta$ or $U_{\theta,T}$),
except that we still need to show \eqref{eq:ess-spec-8}, a task that
we take up in the next section. Hence $R(\ev)$
continues meromorphically to $S_{\pi/4}$ in the
precise sense of the theorem. In particular, $G_o(\ev)$ continues 
meromorphically to $S_{\pi/4}$ as a distribution.

To continue further, we use the following analogue of \eqref{eq:anal-cont-16}:
\begin{equation}\label{eq:anal-cont-24}
\langle f,R(\ev,\theta')g\rangle=\langle U_{\bar\theta}f,R(\ev,\theta+\theta')
U_\theta g\rangle
\end{equation}
for $f,g\in\calA$, $|\im\theta|<\frac{\pi}{4}$. This is 
proved similarly to \eqref{eq:anal-cont-16} using that
$U_\theta U_\theta'=U_{\theta+\theta'}$, and hence
$U_\theta R(\ev,\theta') U_\theta^{-1}=R(\ev,\theta+\theta')$.
Thus, taking $\theta'$ with $\beta<-\im\theta'<\frac{\pi}{4}$,
we can continue $R(\ev,\theta')$ to $S_{-\im\theta'+\frac{\pi}{4}}$
via \eqref{eq:anal-cont-24}.

This gives the extension of $R(\ev)\delta_o$ to $S_{\pi/2}$ as
a distribution. Indeed, we have shown that such an extension exists 
in ${\mathcal D}'(B_T(o))$ for any $T > 0$. The density of $\calA$ 
implies that these extensions are all the same. 
\end{proof}

\section{Analysis of the scaled Laplacian}
We take up, at last, the analysis of the scaled Laplacian $\Delta_{\theta}$.
The key point here is that for each $\theta$ in the appropriate range, 
there is a parametrix construction for $(\Delta_\theta - \ev)^{-1}$ which
is almost the same as that for $R(\ev)$ when $\ev$ is in the resolvent
set. More specifically, the main ingredients in that parametrix construction,
as described in \S 2, are the radial ellipticity of $\Delta$ and the
existence of product-type models at the boundaries of the flat
$\exp(\fraka)$. Both these points have analogues for this scaled 
operator; indeed, we shall discuss the radial ellipticity below,
while the models for $\Delta_\theta$ are obtained by conjugating the 
corresponding models for $\Delta$. For example, the appropriate model 
for $\Delta_\theta$ at $H_\sharp$ is the scaling of $L_\sharp$, 
\begin{equation*}
(L_\sharp)_{\theta}=\frac{1}{4}\,e^{-2\theta}(sD_s)^2+i\frac{1}{2}
\,e^{-\theta}(sD_s)+\frac{1}{3}(\Delta_{\HH^2})_{\theta},
\end{equation*}
which equals the conjugate of $L_\sharp$ by $U_\theta$, and its  
$\SO(2)$-radial part is $(L_\sharp)_{\theta,\mathrm{rad}}$.

The following result is the goal of this parametrix construction. 
\begin{prop}\label{prop:ess-spec}
For $|\im\theta|<\frac{\pi}{2}$, the essential spectrum of 
$\Delta_{\theta,\mathrm{rad}}$ on $L^2(M)^{K_o}$ is given by
\begin{equation*}
\spec_{\mathrm{ess}}(\Delta_{\theta,\mathrm{rad}})=\ev_0+e^{-2i\im\theta}[0,+\infty).
\end{equation*}
\end{prop}

For technical reasons, it is convenient to work with a unitarily equivalent 
problem. Let $\rho$ be a $K_o$-invariant total boundary defining function.
Then  multiplication by $\rho^w$ gives a unitary map
\begin{equation*}
\rho^w:L^2(M,\rho^{2\re w}\,dg)^{K_o}\to L^2(M,dg)^{K_o},
\end{equation*}
hence with
\begin{equation*}
\Deltat_{\theta,\mathrm{rad}}=\rho^{-w}\Delta_{\theta,\mathrm{rad}}\rho^w,\ w=e^\theta,
\end{equation*}
then we have $(\Delta_{\theta,\mathrm{rad}}-\ev)^{-1}=\rho^w(\Deltat_{\theta,\mathrm{rad}}-\ev)^{-1}
\rho^{-w}$. Therefore it suffices to study the analytic continuation of the 
resolvent of $\Deltat_{\theta,\mathrm{rad}}$ instead. The proposition then follows from 
the following lemma. 

\begin{lemma}
For $|\im\theta|<\frac{\pi}{2}$, the essential spectrum of 
$\Deltat_{\theta,\mathrm{rad}}$ on $L^2(M)^{K_o}$ is given by
\begin{equation*}
\spec_{\mathrm{ess}}(\Deltat_{\theta})=\ev_0+e^{-2i\im\theta}[0,+\infty).
\end{equation*}
\end{lemma}

\begin{proof}
Identify $K_o$-invariant functions on $M$ with $S_3$-invariant functions
on $\fraka$, so that $\Deltat_{\theta}$ becomes a three-body
type perturbation of the scaled Euclidean Laplacian on $\fraka$. 
Using language consonant with that interpretation, the subsystem Hamiltonians 
are simply the radial part of the scaled Laplacian 
$(\Delta_{\HH^2})_{\theta}$, conjugated by $\rho_2^w$ (where $\rho_2$ 
is an $\SO(2)$-invariant boundary defining function on $\HH^2$). 
The model for $\Deltat_\theta$ coming from this subsystem is the
$\SO(2)$-radial part of
\begin{equation*}
(L_0)_\theta=\frac{1}{4}e^{-2\theta}(sD_s)^2+\frac{1}{3}\rho_2^{-w}
(\Delta_{\HH^2})_{\theta}\rho_2^{w}+
\frac{1}{4}. 
\end{equation*}
Here $s \in (0,+\infty)$ can be thought of as a variable on 
the wall of a Weyl chamber. For example, at the wall $\lambda_1=\lambda_2$
we take $s=\lambda_3^{-3/2}$. Replacing $s$ by $t=-\log s$ and using the
polar variable $r$ on $\HH^2$, then $(L_0)_{\theta,\mathrm{rad}}$ is a perturbation of 
$e^{-2\theta}(\frac{1}{4}D_t^2+\frac{1}{3}D_r^2)+\frac{1}{3}$ as claimed.
In fact, using the notation from the end of \S 2, $(L_0)_\theta$ is a model 
operator for $\Deltat_\theta$ in the sense that
for $\phi_\sharp$ as in Figure~\ref{fig:flat3},
\begin{equation*}
\phi_\sharp(\Deltat_\theta-(L_0)_\theta)\in\rho_\sharp^w\Diffss^1(\olM).
\end{equation*}
Note in particular that $\rho_\sharp^w$ vanishes at the face $\rho_\sharp=0$ 
when $|\im w|<\frac{\pi}{2}$. Furthermore, as in the introduction, 
$\spec_{\mathrm{ess}}(\rho_2^{-w}(\Delta_{\HH^2})_{\theta,\mathrm{rad}}\, 
\rho_2^{w})= \frac{1}{4}+e^{-2i\im\theta}[0,+\infty)$, and hence
\begin{equation}
\spec_{\mathrm{ess}}(L_0)_{\theta,\mathrm{rad}} = \ev_0 +e^{-2i\im\theta}[0,+\infty).  
\label{eq:essspecl0}
\end{equation}

The parametrix construction can now be carried out {\em exactly} as at the end of 
\S 2. To do this we use (\ref{eq:essspecl0}) and the fact that $\Delta_\theta$ 
remains radially elliptic when $|\im \theta | < \pi/2$. This radial ellipticity 
must be verified on three different regions: near $o \in M$, in a neighbourhood of the 
Weyl chamber walls, 
and in the interior of $\fraka^+$, uniformly out to infinity near the corner 
in $\overline{\fraka}$. In the latter region this is clear since here radial 
ellipticity simply corresponds to the symbol ellipticity of the scaled Euclidean 
Laplacian $e^{-2\theta}\Delta_{\RR^2}$ (in the product $b$-calculus).
Next, near $\frakw_\sharp$ say, we need to check that $L_{\sharp,\theta}$
is symbol elliptic in the edge calculus in a neighbourhood of 
$[0,1)_s \times o'$ in $[0,+\infty) \times \HH^2$, and we have checked this
already in \S 1. Finally, at $o \in M$, the (ordinary!) symbol of
$\Delta_\theta$ is just $e^{-2\theta}|\xi|^2$, and hence it too is invertible
even in a neighbourhood of $o$. This means that for any allowable $\theta$ 
we can choose the supports of $\psi_o$ etc.\ small enough so that the
first step of the parametrix construction can be carried out. 
The second stage of the parametrix construction can also be done because
$\ev$ avoids the essential spectrum of the model operators.
The only step that needs checking is that the analogue of
\cite[Corollary~4.8]{Mazzeo-Vasy:Scattering} still holds, which in
turn only depends on the form of the Schwartz kernel as stated in
\cite[Theorem~9.3]{Mazzeo-Vasy:Resolvents}. This theorem is the consequence of
a stationary phase argument, using the form of the Schwartz kernels of
the resolvents of $\frac{1}{4}e^{-2\theta}(sD_s)^2$
and $\frac{1}{3}\rho_2^{-w}(\Delta_{\HH^2})_{\theta,\mathrm{rad}}\rho_2^{w}$.
These have the same form as the kernels of
the resolvents of $\frac{1}{4}(sD_s)^2$
and $\frac{1}{3}\rho_2^{-w}(\Delta_{\HH^2})_{\mathrm{rad}}\rho_2^{w}$,
provided that the shifted spectral parameter $\ev-\ev_0$ is replaced by
$e^{2\theta}(\ev-\ev_0)$.
In this way we obtain 
a parametrix with compact error for $(\Deltat_{\theta,\mathrm{rad}}-\ev)^{-1}$.
This proves that the essential spectrum is $\frac{1}{3}+e^{-2i\im\theta}[0,+\infty)$.
\end{proof}

This completes the proof of the analytic continuation of
the resolvent on $M$. 

\section{Future directions}
We expect that this same approach will give the meromorphic
continuation of the resolvent of the Laplacian on $\SL(n,\RR)/\SO(n,\RR)$
(and other higher rank symmetric spaces);
this continuation will exist on a Riemann surface ramified at points
corresonding to poles of all `subsystem Hamiltonians', i.e.\ to the poles 
of the analytic continuations of the resolvents on $\SL(k)/\SO(k)$, 
$2\leq k<n$. In fact, the model operators for $\SL(n)/\SO(n)$ are expected 
to be sums of operators, each acting on a different factor in a tensor 
product decomposition, with the operators in each factor being either 
the radial part of the Laplacian on a lower rank symmetric space or a
Euclidean-type operator; this generalizes what we have found here,
where the model is a product $\Real^+_s \times \SL(2)/\SO(2)$.
This should lead to the following prescription: for a given $n$, 
take all non-trivial partitions (e.g.\ for $n=4$, the possibilities
are $2+2$ and $3+1$); for each partition add the poles of the analytic 
continuations of $R(\ev)$ on the corresponding lower rank spaces (possibly
shifted by a constant coming from the Euclidean part, as $(sD_s)^2$ above).
No such branch points occur for $\SL(3)/\SO(3)$ in the region
$\arg(\ev-\ev_0)<\frac{\pi}{2}$, except possibly at the bottom of the 
spectrum, since the continuation of $(\Delta_{\HH^2}-\ev)^{-1}$
has no poles up to the same angle. 

It is certainly conceivable that the continuation of $R(\ev)$ on any 
of the symmetric spaces $\SL(n,\RR)/\SO(n,\RR)$ never has poles (hence
no ramification points either), in $\arg(\ev-\ev_0)<\frac{\pi}{2}$;
it might be possible to prove this using that the residues of these putative
poles would correspond to representations of the group which are known 
not to exist. In any case, this question warrants attention, and this 
discussion illuminates the importance of determining the 
existence of poles of the continued resolvent on $\SL(3)/\SO(3)$.

From the point of view of geometric scattering, this continuation
sets the stage for many directions of investigation. For example,
one should investigate the structure of the scattering operator,
the asymptotic distribution of resonances, trace formulae, etc.

\bibliographystyle{plain}
\bibliography{sm}

\end{document}